\documentclass[10pt,a4paper,reqno]{article}
\usepackage{ae,aecompl}% pour une bonne conversion vers pdf.
\usepackage[cm]{aeguill}%  pour une bonne conversion vers pdf.
\usepackage{a4,pstricks,theorem,amsmath,here,latexsym,amssymb}
\usepackage{epsf,psfrag,psfig,epsfig}
\usepackage[latin1]{inputenc}

\newcommand{\Cat}[1]{\mathcal{C}_{#1}} 
\newcommand{\PS}[1]{\mathcal{S}_{#1}} 
\newcommand{\Ms}{M^\bigstar}

\newcommand{\B}[1]{\overline{#1}}
\newcommand{\orient}[1]{\vec{#1}}
\newcommand{\orientT}[1]{\vec{#1}^T}

\newcommand{\ite}{\noindent $\bullet~$}
\newcommand{\iten}{\noindent -$~$}

\newcommand{\dem}{\noindent \textbf{Proof: }} 
\newcommand{\findem}{\vspace{-.6cm} \begin{flushright} $\square~$ \end{flushright} \vspace{.4cm} }
\newcommand{\findembis}[1]{\hspace{#1} \vspace{.1cm} $\square~$} 

\newtheorem{thm}{Theorem}
\newtheorem{prop}[thm]{Proposition}
\newtheorem{lemma}[thm]{Lemma}
\newtheorem{cor}[thm]{Corollary}
\newtheorem{Def}[thm]{Definition}
\newtheorem{Alg}[thm]{Algorithm}

\title{Bijective counting of tree-rooted maps and shuffles of parenthesis systems}
\author{Olivier Bernardi}
\date{}

\begin{document}

\maketitle
\begin{abstract}
The number of  tree-rooted maps, that is, rooted planar maps with a distinguished spanning tree, of size $n$ is $\Cat{n} \Cat{n+1}$ where $\Cat{n}=\frac{1}{n+1}{2n \choose n}$ is the $n^{th}$ Catalan number.  We present a (long awaited) simple bijection which explains this result. We prove that our bijection is isomorphic to a former recursive construction on shuffles of parenthesis systems due to Cori, Dulucq and Viennot.\\ 
\end{abstract}

\section{Introduction}
In the late sixties, Mullin  published an enumerative result concerning planar maps on which a spanning tree is distinguished \cite{Mullin:tree-rooted-maps}. He proved that the number of rooted planar maps with a distinguished spanning tree, or \emph{tree-rooted maps} for short, of size $n$ is $\Cat{n} \Cat{n+1}$ where $\Cat{n}=\frac{1}{n+1}{2n \choose n}$ is the $n^{th}$ Catalan number. This means that tree-rooted maps of size $n$ are in one-to-one correspondence with pairs of plane trees of size $n$ and $n+1$  respectively. 
But although Mullin asked for a bijective explanation of this result, no natural mapping was found  between tree-rooted maps and pairs of trees. Twenty years later,  Cori, Dulucq and Viennot exhibited one such mapping while working on Baxter permutations \cite{Dulucq:shuffle-parenthesis-system}.  More precisely, they established a bijection between pairs of trees and \emph{shuffles of two parenthesis systems}, that is,  words  on the alphabet $a,\B{a},b,\B{b}$, such that the subword consisting of the letters $a,\B{a}$ and the subword consisting of the letters $b,\B{b}$ are parenthesis systems. It is known that tree-rooted maps are in one-to-one correspondence  with shuffles of two parenthesis systems \cite{Mullin:tree-rooted-maps,Walsh:counting-maps-2}, hence the bijection of Cori \emph{et al.} somehow answers Mullin's question. But this answer is quite unsatisfying in the world of maps. Indeed, the bijection of Cori \emph{et al.} is recursively defined on the set of prefixes of shuffles of parenthesis systems and it was not understood how this bijection could be interpreted on maps. We fill this gap by defining a natural, non-recursive, bijection between tree-rooted maps and pairs made of a tree and a non-crossing partition. Then, we show that our construction is isomorphic to the construction  of Cori \emph{et al.} via the encoding of tree-rooted maps by shuffles of parenthesis systems.\\

%Our construction involves certain orientations on maps towhich we shall refer as \emph{tree-orientations}. These orientations are shown to be in one-to-one correspondence with spanning trees. As a by-product  we obtain a bijection between  spanning forests and \emph{feasible out-degree sequences} on planar graphs.\\

Tree-rooted maps, or alternatively shuffles of parenthesis systems, are in one-to-one correspondence with square lattice walks confined in the quarter plane (we explicit this correspondence in the next section). Therefore, our bijection can also be seen as a way of counting these walks. Some years ago, Guy, Krattenthaler and Sagan worked on walks in the plane \cite{Guy:walks} and exhibited a number of nice bijections. However, they advertised the result of Cori \emph{et al.} as being \emph{considerably harder} to prove bijectively. We believe that the encoding in terms of tree-rooted maps makes this result more natural.\\ 
% More generally, the encoding of lattice walks in tree-rooted maps could be the key for many intriguing enumerative results.\\

%Our bijection also allows the random generation of tree-rooted maps in linear time. It was already the case with Cori \emph{et al.}

The outline of this paper is as follows. In Section \ref{section:preliminary}, we recall some definitions and preliminary results on tree-rooted maps. In Section \ref{section:bijection}, we present our bijection between tree-rooted maps of size $n$ and pairs consisting of a tree and a non-crossing partition of size $n$ and $n+1$ respectively. This simple bijection explains why the number of tree-rooted maps of size $n$ is $\Cat{n}\Cat{n+1}$. In Section  \ref{section:equivalence}, we prove that our bijection is isomorphic to the construction of Cori \emph{et al.}  

Our study requires to introduce a large number of mappings; we refer the reader to Figure \ref{fig:equivalence-bijections} which summarizes our notations. 

%This construction leads us to consider some orientations on maps that we call tree-orientations. In Section \ref{section:out-sequences}, we extend our study of tree-orientations and establish a bijection between spanning forests and feasible out-sequences. In Section  \ref{section:equivalence}, we prove that our bijection is isomorphic to the construction of Cori \emph{et al.} via the encoding of tree-rooted maps by parenthesis-shuffles. Our study requires introducing a great number of mappings; we refer the reader to Figure \ref{fig:equivalence-bijections} which summarizes our notations. 

\section{Preliminary results} \label{section:preliminary}
\begin{figure}[ht!]
\begin{center}
\input{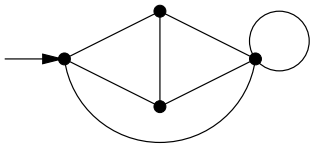} 
\caption{A rooted map.} \label{fig:exemple-carte}
\end{center}
\end{figure}

We begin by some preliminary definitions on planar maps. A \emph{planar map}, or \emph{map} for short, is a two-cell embedding of a connected planar graph into the oriented sphere considered up to orientation preserving homeomorphisms of the sphere. Loops and multiple edges are allowed. A \emph{rooted} map is a map together with a half-edge called the \emph{root}. A rooted map is represented in Figure \ref{fig:exemple-carte}. The vertex (resp. the face) incident to the \emph{root} is called the \emph{root-vertex} (resp. \emph{root-face}). When representing maps in the plane, the root-face is usually taken as the infinite face and the root is represented as an arrow pointing on the root-vertex (see Figure \ref{fig:exemple-carte}). Unless explicitly mentioned, all the maps considered in this paper are rooted. 

\begin{figure}[ht!]
\begin{center}
\hspace{-5cm}\input{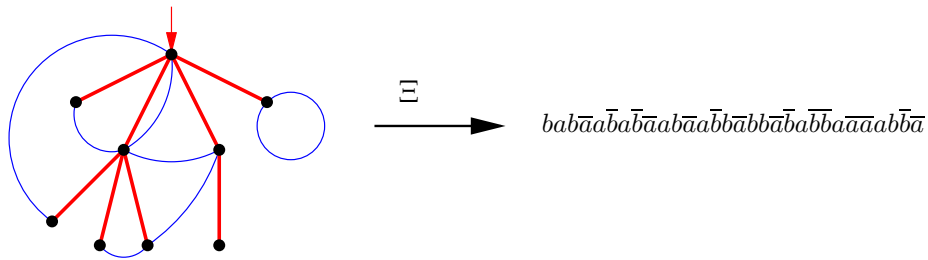} 
\caption{A tree-rooted map and the associated parenthesis-shuffle.} \label{fig:follow-the-border2}
\end{center}
\end{figure}

A \emph{planted plane tree}, or \emph{tree} for short, is a rooted map with a single face. A vertex $v$ is an \emph{ancestor} of another vertex $v'$ in a tree $T$ if $v$ is on the (unique) path in $T$ from $v'$ to the root-vertex of $T$. When $v$ is the first vertex encountered on that path, it is the \emph{father} of $v'$. A \emph{leaf} is a vertex which is not a father. Given a rooted map $M$, a submap of $M$ is a \emph{spanning tree} if it is a tree containing all vertices of $M$. (The spanning tree inherit its root from the map.) We now define the main object of this study, namely tree-rooted maps. A \emph{tree-rooted map} is a rooted map together with  a distinguished spanning tree. 
Tree-rooted maps shall be denoted by symbols like $M_T$ where it is implicitly assumed that $M$ is the underlying map and $T$ the spanning tree.
%If $M$ is a map and $T$ a spanning tree, we denote by $M_T$ the corresponding tree-rooted map. Denoting a tree-rooted map $M_T$ implies implicitly that the underlying map is $M$ and the spanning tree is $T$.  
Graphically, the distinguished spanning tree will be represented by thick lines (see Figure \ref{fig:follow-the-border2}). The \emph{size} of a map, a tree, a tree-rooted map, is the number of edges. \\
\vspace{-.3cm}

\begin{figure}[ht!]
\begin{center}
\hspace{-3cm}\input{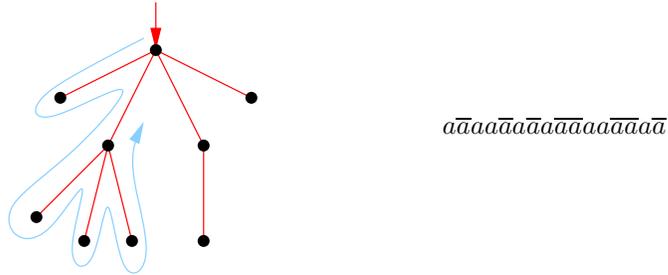} 
\caption{A tree and the associated parenthesis system.} \label{fig:follow-the-border}
\end{center}
\end{figure}

A number of classical bijections on trees are defined by following the border of the tree. Doing the \emph{tour of the tree} means following its border in counterclockwise direction starting and finishing at the root (see Figure \ref{fig:follow-the-border}).
Observe that the tour of the tree induces a linear order, the order of appearance, on the vertex set and on the edge set of the tree. 
For tree-rooted maps,  the tour of the spanning tree $T$ also induces a linear order on half-edges not in $T$ (any of them is encountered once during a tour of $T$). We shall say that a vertex, an edge, a half-edge \emph{precedes} another one \emph{around} $T$. 

%Also, for tree-rooted maps, the tour of the spanning tree $T$ induces an order of appearance on half-edges not in $T$ (any of them is encountered once during a tour of $T$). We shall refer to these orders as the \emph{appearance orders around $T$}.
%More precisely, we say that a vertex $v$ (resp. an edge $e$) of the tree \emph{precedes} another vertex $v'$ (resp. edge $e'$) and we write $v<v'$ (resp. $e<e'$) if it appears before around the tree. Now, consider a tree-rooted map with spanning tree $T$. 
%Making the tour of $T$, any half-edge not in $T$ will be encountered once. So half-edges not in $T$ can be compared according to their order of appearance. %Once again, we say that $h$ \emph{precedes} $h'$ and write $h<h'$ if $h$ appears before $h'$. 
%We shall refer to these orders as the \emph{appearance order} on $T$.

Our constructions lead us to consider \emph{oriented maps}, that is, maps in which all edges are oriented. If an edge $e$ is oriented from  $u$ to $v$, the vertex $u$ is called the \emph{origin} and $v$ the \emph{end}. The half-edge incident to the origin (resp. end) is called the \emph{tail} (resp. \emph{head}). The root of an oriented map will always be considered and represented as a head. \\
\begin{figure}[ht!]
\begin{center}
\input{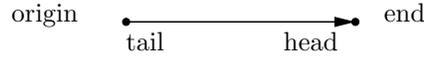}
\caption{Half-edges and endpoints.} \label{fig:head-and-tail}
\end{center}
\end{figure}
\vspace{-.3cm}

%In an oriented map, the \emph{out-degree} of a vertex is the number of tails incident to that vertex. Finally, an \emph{out-degree sequence}, or \emph{out-sequence} for short, for a map is a mapping from its vertices to the set of non-negative integers. An out-sequence $\alpha$ is  \emph{feasible} if there is an orientation such that the out-degree of any vertex $v$ is $\alpha(v)$.\\

We now recall a well-known correspondence between tree-rooted maps and shuffles of two parenthesis systems \cite{Mullin:tree-rooted-maps,Walsh:counting-maps-2}. We derive from it the enumerative result mentioned above: the number of tree-rooted maps of size $n$ (i.e. with $n$ edges) is $\Cat{n} \Cat{n+1}$.  For this purpose, we introduce some notations on words. A word $w$ on a set $A$ (called the alphabet)  is a finite sequence of elements (letters) in $A$. The length of $w$ (that is, the number of letters in $w$) is  denoted $|w|$ and, for $a$ in $A$, the number of occurrences of $a$ in $w$ is denoted  $|w|_a$. A word $w$ on the two-letter alphabet $\{a,\B{a}\}$ is a \emph{parenthesis system} if $|w|_a=|w|_{\B{a}}$ and for all prefixes $w'$, $|w'|_a\geq |w'|_{\B{a}}$. For instance, $aa\B{a}a\B{a}\B{a}$ is a parenthesis system. A \emph{shuffle of two parenthesis systems}, or \emph{parenthesis-shuffle} for short, is a word on the alphabet $\{a,\B{a},b,\B{b}\}$ such that the subword of $w$ consisting of letters in $\{a,\B{a}\}$ and the subword consisting of letters in $\{b,\B{b}\}$ are parenthesis systems. For instance $aba\B{b}\B{a}ba\B{a}\B{b}\B{a}$ is a parenthesis-shuffle. 

\begin{figure}[ht!]
\begin{center}
\input{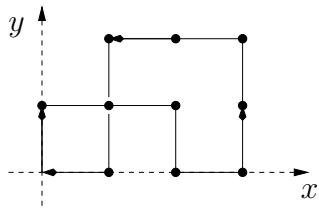}
\caption{A walk in the quarter plane.}\label{fig:square-walks}
\end{center}
\end{figure}

Parenthesis-shuffles can also be seen as walks in the quarter plane. Consider walks made of steps \emph{North, South, East, West,} confined in the quadrant $x\geq 0,~y\geq 0$. The parenthesis-shuffles of size $n$ are in one-to-one correspondence with walks of length $2n$ starting and returning at the origin. This correspondence is obtained by considering each letter $a$ (resp. $\B{a},b,\B{b}$) as a \emph{North} (resp. \emph{South, East, West}) step. For instance, we represented the walk corresponding to $abb\B{a}baa\B{b}\B{b}\B{a}\B{a}\B{b}$ in Figure \ref{fig:square-walks}. The fact that the subword of $w$ consisting of letters in $\{a,\B{a}\}$ (resp. $\{b,\B{b}\}$) is a parenthesis system implies that the walk stays in the half-plane $y\geq 0$ (resp. $x\geq 0$) and returns at $y=0$ (resp. $x=0$).\\

The size of a parenthesis system, a parenthesis-shuffle, is half its length. For instance, the parenthesis-shuffle $aba\B{b}\B{a}ba\B{a}\B{b}\B{a}$ has size 5.  It is well known that the number of parenthesis systems of size $n$ is the $n^{th}$ Catalan number $\Cat{n}=\frac{1}{n+1}{2n \choose n}$. From this, a simple calculation proves that the number of parenthesis-shuffles of size $n$ is $\PS{n}=\Cat{n} \Cat{n+1}$. Indeed, there are ${2n \choose 2k}$ ways to shuffle a parenthesis system of size $k$ (on  $\{a,\B{a}\}$)  with a parenthesis system of size $n-k$ (on  $\{b,\B{b}\}$). And summing on $k$ gives the result:
$$
\hspace{0cm}\begin{array}{ll}
\displaystyle \PS{n} &=~\displaystyle \sum_{k=0}^n {2n \choose 2k} \Cat{k}\Cat{n-k}~=~\frac{(2n)!}{(n+1)!^2}\sum_{k=0}^n {n+1 \choose k} {n+1 \choose n-k}\vspace{.2cm} \\
\displaystyle &=~ \displaystyle \frac{(2n)!}{(n+1)!^2}{2n+2 \choose n} ~=~ \Cat{n}\Cat{n+1}.
\end{array} 
$$
Note, however, that this calculation involves the Chu-Vandermonde identity.

It remains to show that tree-rooted maps of size $n$ are in one-to-one correspondence with parenthesis-shuffles of size $n$. We first recall a very classical bijection between trees and parenthesis systems. This correspondence is obtained by making the tour of the tree. Doing so and writing $a$ the first time we follow an edge and $\B{a}$  the second time we follow that edge (in the opposite direction) we obtain a parenthesis system. This parenthesis system is indicated for the tree of Figure \ref{fig:follow-the-border}. Conversely, any parenthesis system can be seen as a code for constructing a tree. 

%In the sequel, we shall simply say that we \emph{follow the border of a tree} to mean that we follow it  in counterclockwise direction starting and finishing at the root. 
\noindent Now, consider  a tree-rooted map. During the tour of the spanning tree we cross edges of the map that are not in the spanning tree. In fact, each edge not in the spanning tree will be crossed twice (once at each half-edge). Hence, making the tour of the spanning tree and  writing $a$ the first time we follow an edge of the tree, $\B{a}$  the second time,  $b$ the first time we cross an edge not in the tree and $\B{b}$  the second time, we obtain a parenthesis-shuffle. We shall denote by $\Xi$ this mapping from tree-rooted maps to parenthesis-shuffles. We applied the mapping $\Xi$ to the tree-rooted map of Figure \ref{fig:follow-the-border2}.% (the spanning tree is indicated by thick lines).

\noindent The reverse mapping can be described as follows: given a parenthesis-shuffle $w$ we first create the tree corresponding to the subword of $w$ consisting of letters $a,\B{a}$ (this will give the spanning tree) then we glue to this tree a head for each letter $b$ and a tail for each letter $\Bar{b}$. There is only one way to connect heads to tails so that the result is a planar map (that is, no edges intersect). Note that, if the map $M$ has  size $n$, the corresponding parenthesis-shuffle $w$ has size $n$ since $|w|_{a}$ is the number of edges in the tree and  $|w|_{b}$ is the number of edges not in the tree.\\
This encoding due to Walsh and Lehman \cite{Walsh:counting-maps-2} establishes a one-to-one correspondence between  tree-rooted maps of size $n$ and parenthesis-shuffles of size $n$. Hence, there are $\Cat{n}\Cat{n+1}$ tree-rooted maps of size $n$. \\

Such an elegant enumerative result is intriguing for combinatorists since Catalan numbers have very nice combinatorial interpretations. We have just seen that these numbers count parenthesis systems and trees. In fact, Catalan numbers appear in many other contexts (see for instance Ex. 6.19 of \cite{Stanley:volume2} where 66 combinatorial interpretations are listed). We now give another classical combinatorial interpretation of Catalan numbers, namely \emph{non-crossing partitions}. A non-crossing partition is an equivalence relation $\sim$ on a linearly ordered set $S$ such that no elements $a<b<c<d$ of $S$ satisfy $a \sim c$, $b \sim d$ and  $a \nsim b$. The equivalence classes of non-crossing partitions are called \emph{parts}. Non-crossing partitions have been extensively studied (see \cite{Simion:non-crossing-partition} and references therein). 
%The size of a non-crossing partition is the size of the set $S$. The equivalence classes of non-crossing partitions are called \emph{parts}. Non-crossing partitions have been extensively studied (see \cite{Simion:non-crossing-partition} and references therein). 

\begin{figure}[ht!]
\begin{center}
\input{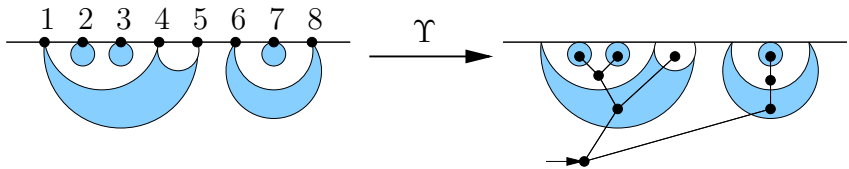} 
\caption{A non-crossing partition and the associated tree.}\label{fig:partition-tree}
\end{center}
\end{figure}

Non-crossing partitions can be represented as cell decompositions of the half-plane.
If the set $S$ is $\{s_1,\ldots,s_n\}$ with $s_1<s_2<\cdots<s_n$, we associate with $s_i$ the vertex of coordinates $(i,0)$ and  with each part we associate a connected region of the lower half-plane $y\leq 0$ incident to the vertices of that part. The existence  of a cell decomposition with no intersection between cells is precisely the definition of non-crossing partitions.  A non-crossing partition of size 8 is represented in Figure \ref{fig:partition-tree}. The only non-trivial parts of this non-crossing partition are $\{1,4,5\}$ and $\{6,8\}$. 

Non-crossing partitions of size $n$ (i.e. on a set of size $n$) are in one-to-one correspondence with trees of size $n$. One way of seeing this is to draw the dual of the cell-representation of the partition, that is, to draw a vertex in each part and each \emph{anti-part} (connected cells complementary to parts in the half-plane decomposition) and connect vertices corresponding to adjacent cells by an edge. The root is chosen in the infinite cell as indicated in Figure \ref{fig:partition-tree}. In the sequel, this mapping between non-crossing partitions and trees is denoted $\Upsilon$. It is a bijection between  non-crossing partitions of size $n$ and trees of size $n$. It proves that the number of non-crossing partitions of size $n$ is $\Cat{n}$.\\

\section{Bijective decomposition of tree-rooted maps } \label{section:bijection}
We begin with the presentation of our bijection between tree-rooted maps and  pairs consisting of a tree and a non-crossing partition. This bijection has two steps: first we orient the edges of the map and then we disconnect properly the vertices.\\

\noindent \textbf{Map orientation:} Let $M_T$ be a tree-rooted map. We denote by  $\orientT{M}$ the oriented map obtained by orienting the edges of $M$ according to the following rules:\\
%Let $M$ be a rooted map and $T$ a distinguished spanning tree. 
%We associate an oriented map $\orientT{M}$  with the tree-rooted map  $M_T$ by orienting its edges as follows: 
\ite edges in the tree $T$ are oriented from the root to the leaves,\\
\ite edges not in the tree $T$ are oriented in such a way that their head precedes their tail around $T$.\\
As always in this paper, the root is considered as a head. 

\begin{figure}[ht!]
\begin{center}
\input{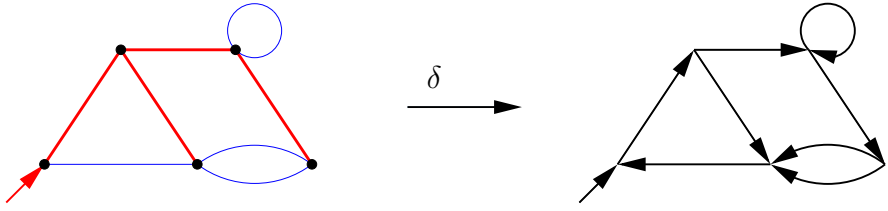} 
\caption{A tree-rooted map $M_T$ and the corresponding oriented map $\orientT{M}$.} \label{fig:orientation}
\end{center}
\end{figure}

In the sequel, the mapping $M_T \mapsto \orientT{M}$  is denoted $\delta$. We applied this mapping to the tree-rooted map of Figure \ref{fig:orientation}. 
Note that any vertex of $\orientT{M}$ is incident to at least one head (since the spanning tree is oriented from the root to the leaves).\\

\begin{figure}[ht!]
\begin{center}
\input{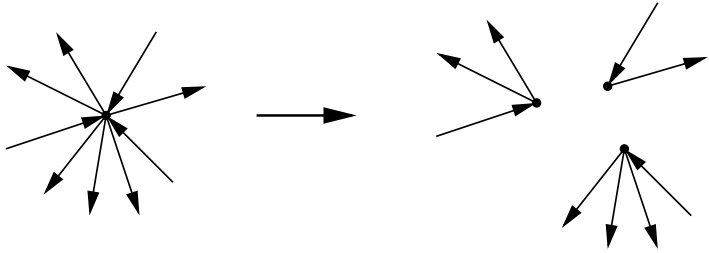} 
\caption{Local rule for suppressing the adjacency relations.}\label{fig:vertex-explosion}
\end{center}
\end{figure}

\noindent \textbf{Vertex explosion:} We replace each vertex $v$ of the oriented map $\orientT{M}$ by as many vertices as heads incident to $v$ and we suppress some adjacency relations between half-edges incident to $v$ according to the rule represented in Figure \ref{fig:vertex-explosion}. That is, each tail $t$ becomes adjacent to exactly one head which is the first head encountered in counterclockwise direction around  $v$ starting from $t$.\\

%Given the oriented map $\orientT{M}$, we suppress some adjacency relations between edges. More precisely,  we replace each vertex $v$ of $\orientT{M}$ by as many vertices as heads incident to $v$ and we modify the adjacency relations between half-edges incident to $v$ according to the rule represented in Figure \ref{fig:vertex-explosion}. That is, we suppress some adjacency relations between half-edges incident to $v$ in such a way, each tail $t$ becomes adjacent to exactly one head which is the first head encountered in counterclockwise direction around  $v$ starting from $t$.\\

We shall prove  (Lemma \ref{thm:explosion-produce-tree}) that this suppression of some adjacency relations in $\orientT{M}$ produces a tree denoted $\varphi_0(\orientT{M})$. Observe that this  tree  has the same number of edges, say $n$, as the original map $M$. Hence, its vertex set  $S$ has size $n+1$. This set is linearly ordered by the order of appearance around the tree  $\varphi_0(\orientT{M})$.
We define an equivalence relation $\varphi_1(\orientT{M})$ on $S$: two vertices are equivalent if they come from the same vertex of $\orientT{M}$.  We will prove (Lemma \ref{thm:explosion-produce-partition}) that the equivalence relation  $\varphi_1(\orientT{M})$ is a non-crossing partition on the set $S$. The mapping $\orientT{M}  \mapsto  (\varphi_0(\orientT{M}),\varphi_1(\orientT{M}))$ is called  the \emph{vertex explosion process} and is denoted $\varphi$.\\

Therefore, with any tree-rooted map $M_T$ of size $n$  we associate  a tree $\varphi_0(\orientT{M})$ of size $n$ and a non-crossing partition $\varphi_1(\orientT{M})$ of size $n+1$.  The following theorem states that this correspondence is one-to-one.

\begin{thm} \label{thm:phi-bijective}
Let $\Phi$ be the mapping associating the ordered pair $(\varphi_0(\orientT{M}),\varphi_1(\orientT{M}))$ with the tree-rooted map $M_T$. This mapping is a bijection  between the set of tree-rooted maps of size $n$ and the Cartesian product of the set of trees of size $n$ and the set of  non-crossing partitions of size $n+1$.\\
It follows that the number of tree-rooted maps of size $n$ is  $\Cat{n}\Cat{n+1}$.
\end{thm} 

%\begin{thm} \label{thm:phi-bijective}
%Let $\Phi$ be the mapping associating the ordered pair $(\varphi_0(\orientT{M}),\varphi_1(\orientT{M}))$ with a tree-rooted map $M_T$. The mapping $\Phi$ is a bijection  between the set $\TM{n}$ of tree-rooted maps of size $n$ and the Cartesian product $\Tr{n}\times \Pa{n+1}$ of the set $\Tr{n}$ of trees of size $n$ and the set $\Pa{n+1}$ of  non-crossing partitions of size $n+1$.\\
%It follows that the number of tree-rooted maps of size $n$ is  $|\TM{n}|= \Cat{n}\Cat{n+1}$.
%\end{thm} 

\begin{figure}[ht!]
\begin{center}
\input{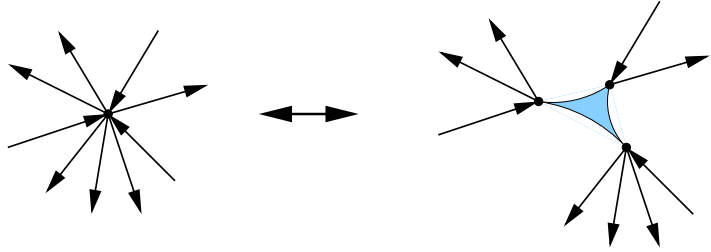} 
\caption{The vertex explosion process and a part of the non-crossing partition.}\label{fig:vertex-explosion2}
\end{center}
\end{figure}

\begin{figure}[ht!]
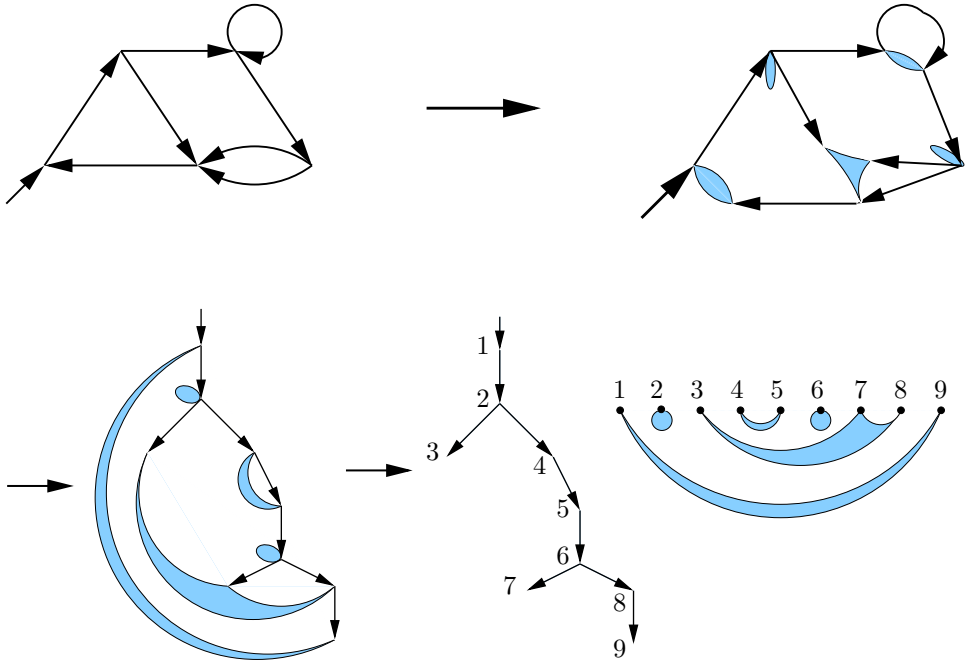

\begin{center}
\hspace{-.3cm} \input{exemple2bis.pstex_t} \hspace{-.3cm} \input{exemple3bisbis.pstex_t}
\caption{The vertex explosion process $\varphi$.} \label{fig:exemple-phi}
\end{center}
\end{figure}

Graphically, the bijection $\Phi$ is best represented by keeping track of the underlying non-crossing partition during the vertex explosion process. This is done by creating for each vertex of $M$ a connected cell representing the corresponding part of the non-crossing partition. The graphical representation of the vertex explosion process $\varphi$ becomes as indicated in Figure \ref{fig:vertex-explosion2}. For instance, we applied the mapping $\varphi$ to the oriented map of Figure \ref{fig:exemple-phi}.  \\

\begin{figure}[ht!]
\begin{center}
\hspace{0cm} \input{exemple1bis.pstex_t}
\end{center}
\end{figure}

The rest of this section is devoted to the proof of Theorem \ref{thm:phi-bijective}. We first give a characterization of the set of oriented maps, called \emph{tree-oriented maps}, associated to tree-rooted maps by the mapping $\delta$. We also define the reverse mapping $\gamma$.  Then we prove that the vertex explosion process $\varphi$ is a bijection between tree-oriented maps (of size $n$) and pairs made of a tree and a non-crossing partition (of size $n$ and $n+1$ respectively).\\

\subsection{Tree-rooted maps and tree-oriented maps} 
In this subsection, we consider certain orientations of maps called \emph{tree-orientations} (Definition \ref{def:tree-orientation}).
We prove that the mapping $\delta: M_T \mapsto \orientT{M}$ restricted to any given map $M$ induces a bijection between spanning trees and tree-orientations of $M$. The key property explaining why the mapping $\delta$ is injective is that during a tour of a spanning tree $T$, the tails of edges in $T$ are encountered before their heads whereas it is the contrary for the edges not in $T$. Using this property we will define a procedure $\gamma$ for recovering spanning trees from tree-orientations of $M$ (Definition \ref{def:recover-tree}).  We will prove that $\delta$ and $\gamma$ are reverse mappings that establish a one-to-one correspondence between tree-rooted maps and tree-oriented maps (Proposition \ref{thm:tree-coloration-orientation}).\\

\begin{figure}[ht]
\begin{center}
\hspace{0cm} \input{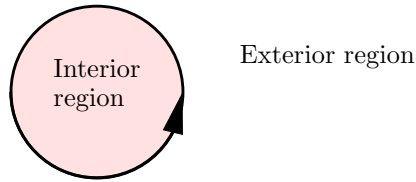}
\caption{Interior and exterior regions of a directed cycle.}\label{fig:interior-region}
\end{center}
\end{figure}

We begin with some definitions concerning cycles and paths in oriented maps. 
A simple cycle (resp. simple path) is \emph{directed} if all its edges are oriented consistently.  A simple cycle defines two regions of the sphere. The \emph{interior region} (resp. \emph{exterior region}) of a directed cycle is the region situated at its left (resp. right) as indicated in  Figure \ref{fig:interior-region}. We call \emph{positive cycle} a directed cycle having the root in its exterior region. Graphically, positive cycles appear as counterclockwise directed cycles when the map is projected on the plane with the root in the infinite face.

\begin{Def}\label{def:tree-orientation}
%An orientation of a map is a 
A \emph{tree-orientation} of a map is an orientation without positive cycle such that any vertex can be reached from the root by a directed path. A \emph{tree-oriented map} is a map with a tree-orientation.
\end{Def}

We will prove that the images of tree-rooted maps by the mapping $\delta$ are tree-oriented maps. More precisely, we have the following proposition.
\begin{prop}\label{thm:tree-coloration-orientation}
For any given map $M$, the mapping $\delta: M_T \mapsto \orientT{M}$ induces a bijection between spanning trees and tree-orientations of $M$.
\end{prop}

We first prove the following lemma.

\begin{lemma}\label{thm:lemme-delta}
For all tree-rooted map $M_T$, the map $\orientT{M}$ is tree-oriented.
\end{lemma}

\begin{figure}[ht!]
\begin{center}
\input{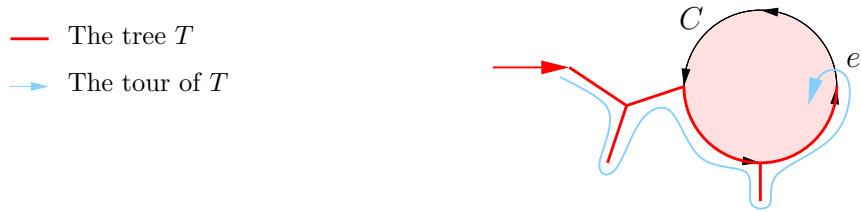}
\caption{Entering the cycle $C$.} \label{fig:tree-oriented}
\end{center}
\end{figure}

\dem For any vertex $v$, there is a  path in $T$ from the root to $v$. This path is oriented from the root to $v$ in  $\orientT{M}$. It remains to prove that there is no positive cycle. 
Suppose the contrary and consider a positive cycle $C$.  By definition, the root is in the exterior region of $C$. Since $C$ is a cycle there are edges of $C$ which are not in $T$. Consider the first such edge $e$ encountered during the tour of $T$. When we first cross $e$ we enter for the first time the interior region of $C$. Given the orientation of $C$, the half-edge of $e$ that we first cross is its tail (see Figure \ref{fig:tree-oriented}). But, by definition of $\orientT{M}$, the half-edge of $e$ that we first cross should be its head. This gives a contradiction.
\findem

We now define a procedure $\gamma$ constructing a spanning tree $T$ on a tree-oriented map $\orient{M}$. %For the need of our definition an edge is either \emph{active} or \emph{inactive}.

\begin{Alg}\label{def:recover-tree}
\vspace{.1cm} ~\\
\framebox{\parbox{12.5cm}{
\textbf{\hspace{.1cm} Procedure $\gamma$:}
\begin{enumerate}
\item At the beginning, the submap $T$ is reduced to the root and the root-vertex.  
\item We make the tour of $T$ (starting from the root) and apply the following rule. \\ 
When the tail of an edge $e$ is encountered  and its head has not been encountered yet, we add $e$ to $T$ (together with its end).\\
Then we continue the tour of $T$, that is, if $e$ is in $T$ we follow its border, otherwise we cross $e$.
\item We stop when arriving at the root and return the submap $T$.
\end{enumerate}}}
\end{Alg}

%Note that an edge is active if and only if none of its half-edge has been encountered. 
We prove the correction of the procedure  $\gamma$. 
\begin{lemma} \label{thm:lemme-gamma}
The mapping $\gamma$ is well defined (terminates) on tree-oriented maps and returns a spanning tree. 
\end{lemma}

\begin{figure}[ht!]
\begin{center}
\input{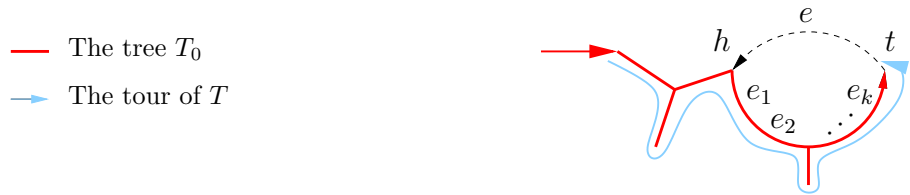} 
\caption{The submap $T$ remains a tree.}\label{fig:dem-tree}\vspace{-0cm}
\end{center}
\end{figure}

\dem \\
\ite \emph{At any stage of the procedure, the submap $T$ is a tree.} \\
Suppose not, and consider the first time an edge $e$ creating a cycle is added to $T$. We denote by $T_0$ the tree $T$ just before that time. The edge $e$ is added to $T_0$ when its tail $t$ is encountered. At that time, its head $h$ has not been encountered but is incident to $T_0$ (since adding $e$ creates a cycle). We know that, when $e$ is added, the border of $T_0$ from the root to $t$  has been followed but not the border of  $T_0$ from $t$ to the root. Moreover, the head $h$ lies after $t$ around $T_0$ (since $h$ has not been encountered yet).  Observe that the right border of any edge of $T_0$ has been followed (just after this edge was added to $T_0$). Thus, the border of  $T_0$ from $t$ to $h$ is made of the left borders of some edges $e_1,e_2,\ldots,e_k$. Hence, these edges  form a directed path from $h$ to $t$ and  $e,e_1,e_2,\ldots,e_k$ form a directed cycle $C$. Since $h$ lies after $t$ around $T_0$, the root is in the exterior region of $C$ (see Figure \ref{fig:dem-tree}). Therefore, the cycle $C$ is positive which is impossible.\\

\ite   \emph{The procedure $\gamma$ terminates.} \\
The set $T$ remains a tree connected to the root. Hence, it is impossible to follow the same border of the same edge twice without encountering the root.\\
\ite \emph{At the end of the procedure $\gamma$, the tree $T$ is spanning.}\\ 
At the end of the procedure, the whole border of $T$ has been followed. Hence, any half-edge incident to $T$ has been encountered. 
Now, suppose that a vertex $v$ is not in $T$ and consider a directed path from the root to $v$. (This path exists by definition of tree-orientations.) There is an edge of this path with its origin in $T$ and its end out of $T$. Therefore, its tail is incident to $T$ but not its head. Thus, it should have been added to $T$ (with its end) when its tail was encountered. This is a contradiction.$~\square$
%\findem

We continue the proof of Proposition \ref{thm:tree-coloration-orientation}. We proved that the mapping $\delta$ associates a tree-orientation of a map to any spanning tree of that map (Lemma \ref{thm:lemme-delta}). We proved that the mapping $\gamma$ associates a spanning tree of a map to any tree-orientation of that map (Lemma \ref{thm:lemme-gamma}). It remains to prove that $\delta \circ \gamma$ and $\gamma \circ \delta$ are identity mappings. \\
%We say that a map is  \emph{canonically oriented}  with respect to a spanning tree $T$ if edges of $T$ are oriented from the root to the leaves, edges not in $T$ are oriented in such a way that their head precedes their tail. By definition the map $\orientT{M}$ is canonically oriented with respect to $T$. We prove a last lemma.

\begin{lemma} \label{thm:tree+orientation}
Let $\orient{M}$ be a tree-oriented map and $T$ be the spanning tree constructed by the procedure $\gamma$. The edges in $T$ are oriented from the root to the leaves and the edges not in $T$ are oriented in such a way that their head precedes their tail around~$T$.
%The tree-oriented map $\orient{M}$ is canonically oriented with respect to the tree $T$ constructed on $\orient{M}$ by the procedure $\gamma$.
\end{lemma}

\dem \\
\ite \emph{Edges in $T$ are oriented from the root to the leaves.} 
An edge $e$ is added to $T$ when its tail is encountered. At that time the end of $e$ is not in $T$ or adding $e$ would create a cycle. The property follows by induction.\\
\ite \emph{Edges not in $T$ are oriented in such a way that their head precedes their tail around $T$.} 
If an edge breaks this rule it should have been added to $T$ when its tail was encountered.
\findem

\begin{cor}
The mapping $\delta \circ \gamma$ is the identity mapping on tree-oriented maps.
\end{cor}

\dem Let $\orient{M}$ be a tree-oriented map and $T$ be the tree constructed by the procedure $\gamma$. By Lemma \ref{thm:tree+orientation},  the edges in $T$ are oriented from the root to the leaves and the edges not in $T$ are oriented in such a way that their head precedes their tail around $T$. By definition of $\delta$, this is also the case in $\delta\circ \gamma(\orient{M})$. Thus,  $\delta \circ \gamma$ is the identity mapping on tree-oriented maps. 
\findem

%\dem Let $\orient{M}$ be a tree-oriented map and $T$ the tree constructed by procedure $\gamma$. From Lemma \ref{thm:tree+orientation}  the map $\orient{M}$ is canonically oriented with respect to $T$. By definition, the  map $\orientT{M}=\delta\circ \gamma(\orient{M})$ is canonically oriented with respect to $T$. Therefore,  $\orientT{M}=\orient{M}$. Thus,  $\delta \circ \gamma$ is the identity mapping on tree-oriented maps.
%\findem

\begin{lemma} \label{thm:tree+orientation2}
The mapping $\gamma \circ \delta$ is the identity mapping on tree-rooted maps.
\end{lemma}

\dem Let $M_T$ be a tree-rooted map. Suppose the spanning tree $T'$ constructed by the procedure $\gamma(\delta(M_T))$ differs from $T$. We consider the order of edges induced by the tour of $T$. Let $e$ be the smallest edge in the symmetric difference of $T$ and $T'$. The tours of $T$ and $T'$ must coincide until a half-edge $h$ of $e$ is encountered. We distinguish the head and the tail of $e$ according to its orientation in $\delta(M_T)$. If $e$ is in $T$, its tail is encountered before its head around $T$ (by definition of $\delta(M_T)$). In this case, $h$ is a tail. If $e$ is not in $T'$, its head is encountered before its tail around $T'$ (by Lemma \ref{thm:tree+orientation}). In this case, $h$ is a head. Therefore, $e$ cannot be in $T\setminus T'$. Similarly, $e$ cannot be in $T'\setminus T$ since $e$ being in $T'$ implies that $h$ is a head and $e$ not being in $T$ implies that $h$ is a tail. We obtain a contradiction.
\findem

%\dem Let $M_T$ be a tree-rooted map. Suppose the spanning tree $T'$ constructed by the procedure $\gamma(\delta(M_T))$ differs from $T$. The tours of $T$ and $T'$ must coincide until a certain edge in $T\setminus T' \cup T'\setminus T $. Consider an edge $e$ in $T\setminus T'$.  We distinguish the head and the tail of $e$ according to its orientation in $\delta(M_T)$. Since $e$ is in $T$, its tail is encountered before its head around $T$ (by definition of $\delta(M_T)$). Since $e$ is not in $T'$, its head is encountered before its tail around $T'$ (by Lemma \ref{thm:tree+orientation}). Therefore $e$ is not the first edge at which the tours of $T$ and $T'$ differ. Similarly, no edge of $T'\setminus T$ can be the first edge at which the tours of $T$ and $T'$ differ. We have a contradiction.
%\findem

This completes the proof of Proposition \ref{thm:tree-coloration-orientation}: tree-oriented maps are in one-to-one correspondence with tree-rooted maps. \findembis{6cm}
%\findem

\subsection{The  vertex explosion process on tree-oriented maps} \label{subsection:vertex-explosion}
This subsection is devoted to the proof of the following proposition.
\begin{prop}\label{thm:varphi-bijective} The mapping $\varphi: \orient{M} \mapsto  (\varphi_0(\orient{M}),\varphi_1(\orient{M}))$ is a bijection between tree-oriented maps of size $n$ and ordered pairs consisting of a tree of size $n$ and a non-crossing partition of size $n+1$. 
\end{prop}

We start with a  lemma concerning  the mapping $\varphi_0$. 
\begin{lemma} \label{thm:explosion-produce-tree} 
The image of any  tree-oriented map $\orient{M}$ by $\varphi_0$ is a tree (oriented from the root to the leaves).
%Let $\varphi_0$ be the mapping defined on tree-oriented maps by the vertex explosion process illustrated in Figure \ref{fig:vertex-explosion}.  The mapping $\varphi_0$  returns a tree (oriented from the root to the leaves).
\end{lemma}

\dem 
Let $\orient{M}$ be a tree-oriented map. %It is not clear, at first sight, that $\varphi_0(\orient{M})$ is connected. 
Any vertex is incident to at least one head (there is a directed path from the root to any vertex), hence the mapping $\varphi_0$ is well defined. The image $\varphi_0(\orient{M})$ has the same number of edges, say $n$, as $\orient{M}$. The map $\orient{M}$ has $n+1$ heads (one per edge plus one for the root). Since any vertex in  $\varphi_0(\orient{M})$ is incident to exactly one head, the image $\varphi_0(\orient{M})$ has $n+1$ vertices. Thus, it is sufficient to prove that $\varphi_0(\orient{M})$ has no cycle (it will imply the connectivity).\\
Suppose $\varphi_0(\orient{M})$ contains a simple cycle $C$. Since any vertex in  $C$ is incident to exactly one head, the edges of $C$ are oriented consistently. We identify the edges of $\orient{M}$ and the edges of $\varphi_0(\orient{M})$. The edges of $C$ form a cycle  in $\orient{M}$ but this cycle might not be simple. We consider a directed path $P$ in $\orient{M}$ from the root to a vertex $v$ (of $\orient{M}$) incident with an edge of $C$. We suppose (without loss of generality) that $v$ is the only vertex of $P$ incident with an edge of $C$. Let $h$ be the head in $P$ incident with $v$ and $t'$ be the first tail in $C$ following $h$ in counterclockwise direction around $v$. We can construct a directed simple cycle $C'$ (in $\orient{M}$) made of edges in $C$ and containing $t'$ (see Figure \ref{fig:cycle-explosion}). Let $h'$ be the head of $C'$ incident with $v$. Since $C'$ is a directed cycle of the tree-oriented map $\orient{M}$, it contains the root in its interior region. Since $v$ is the only vertex of $P$ incident with an edge in $C'$, the head $h$ is in the interior region of $C'$. Therefore, in counterclockwise direction around $v$ we have $h,~h'$ and $t'$ (and possibly some other half-edges). We consider the tail $t$ following $h$ in the cycle $C$ (considered as a directed simple cycle of $\varphi_0(\orient{M})$). By the choice of $t'$ we know that $t$ is between $t'$ and $h$ in counterclockwise direction around $v$ ($t$ and $t'$ may be distinct or not). Hence, in counterclockwise direction around $v$ we have $h,~h'$ and $t$. Hence,  $h'$ is not the first head encountered in counterclockwise direction around $v$ starting from $t$. Therefore, by definition of the vertex explosion process, $h'$ and $t$ are not adjacent in  $\varphi_0(\orient{M})$. We reach a contradiction.
\findem

\begin{figure}[ht!]
\begin{center}
\input{cycle-explosion2.pstex_t} 
\caption{The cycle $C'$ in $\orient{M}$.} \label{fig:cycle-explosion}
\end{center}
\end{figure}

\begin{figure}[ht!]
\begin{center}
\input{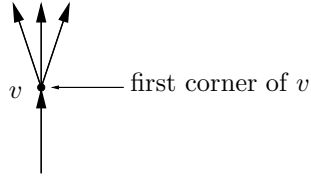}
\caption{The first corner of a vertex.} \label{fig:first-corner}
\end{center}
\end{figure}

We now study the properties of the mapping $\varphi_1$. Two consecutive half-edges around a vertex define a \emph{corner}. A vertex has as many corners as incident half-edges.  Let $T$ be a tree and $v$ be a vertex of $T$. The \emph{first corner} of the vertex $v$ is the first corner of $v$ encountered around $T$. If the tree is oriented from the root to the leaves, the first corner of $v$ is at the right of the head incident to $v$ as shown in Figure \ref{fig:first-corner}.

\noindent
%More precisely, we say that a vertex $v$ (resp. an edge $e$) of the tree \emph{precedes} another vertex $v'$ (resp. edge $e'$) and we write $v<v'$ (resp. $e<e'$) if it appears before around the tree. Now, consider a tree-rooted map with spanning tree $T$. 
%Making the tour of $T$, any half-edge not in $T$ will be encountered once. So half-edges not in $T$ can be compared according to their order of appearance. %Once again, we say that $h$ \emph{precedes} $h'$ and write $h<h'$ if $h$ appears before $h'$. 
%We shall refer to these orders as the \emph{appearance order} on $T$.

We compare the vertices of the tree $\varphi_0(\orient{M})$ according to their order of appearance around this tree. We write  $u<v$  if $u$ precedes $v$ (i.e. the first corner of $u$ precedes the first corner of $v$) around the tree. 
%that is, $v<v'$ if the  first corner of $v$ appears before the first corner of $v'$ around the tree. 

%We now state a lemma concerning  the equivalence relation $\varphi_1(orient{M})$.
\begin{lemma} \label{thm:explosion-produce-partition}
For any tree-oriented map  $\orient{M}$, the equivalence relation $\varphi_1(\orient{M})$ on the set of vertices of the tree $\varphi_0(\orient{M})$ ordered by their order of appearance around this tree is a non-crossing partition.
\end{lemma}

\begin{figure}[ht!]
\begin{center}
\input{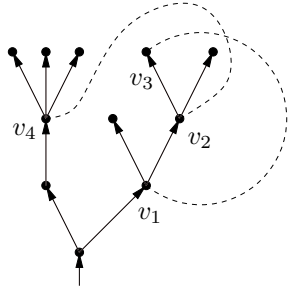}
\caption{The two paths  intersect.} \label{fig:intersect}
\end{center}
\end{figure}

\dem The proof relies on  the graphical representation of the equivalence relation  $\sim =\varphi_1(\orient{M})$ given by Figure \ref{fig:vertex-explosion2}. During the vertex explosion process, we associate a connected cell $C_v$ with each vertex $v$ of $\orient{M}$, that is, with each equivalence class of the relation  $\sim$. The cell $C_v$ can be chosen to be incident only with the first corners of the vertices in its class but not otherwise incident with the tree. Moreover the cells can be chosen so that they do not intersect.\\
Suppose $v_1<v_2<v_3<v_4$, $v_1\sim v_3$ and $v_2\sim v_4$. One can draw a path from the first corner of $v_1$ to the first corner of $v_3$ staying in a cell $C$ and a path from the first corner of $v_2$ to the first corner of $v_4$ staying in a cell $C'$. It is clear that these two paths  intersect (see Figure \ref{fig:intersect}). Thus $C=C'$ and $v_1\sim v_2$. 
\findem

We have proved that the application $\varphi: \orient{M} \mapsto (\varphi_0(\orient{M}),\varphi_1(\orient{M}))$ associates a tree of size $n$ and a non-crossing partition of size $n+1$ with any tree-oriented map of size $n$.  Conversely, we define  the mapping $\psi$.

\begin{Def} \label{def:psi}
Let $T$ be a tree of size $n$ and $\sim$ be a non-crossing partition on a linearly ordered set $S$ of size $n+1$. We identify $S$ with the set of vertices of $T$ ordered by the order of appearance around $T$. We construct the  oriented map $\psi(T,\sim)$ as follows. 
First we orient the tree $T$ from the root to the leaves.  With each part $\{v_1,v_2,\ldots,v_k\}$ of the partition, we associate a simply connected cell incident to the first corner of $v_i, ~i=1\ldots k$ but not otherwise incident with $T$. Since $\sim$ is a non-crossing partition, these cells can be chosen without intersections.  Then we contract each cell into a vertex in such a way  no edges of $T$ intersect. 
%The oriented map obtained is denoted $\psi(T,\sim)$.
\end{Def}

We first prove the following lemma.
\begin{lemma}\label{thm:psi-produces-tree-oriented}
For any tree $T$ of size $n$ and any non-crossing partition $\sim$ of size $n+1$, the oriented map $\psi(T,\sim)$ is tree-oriented.
\end{lemma}

\begin{figure}[ht!]
\begin{center}
\input{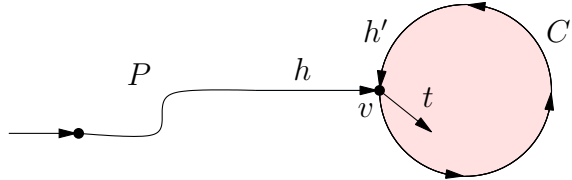}
\caption{The map $\orient{M}=\psi(T,\sim)$ has no positive cycle.} \label{fig:dem-merging}
\end{center}
\end{figure}

\dem
Every vertex of $\orient{M}=\psi(T,\sim)$ is connected to the root by a directed path (since it is the case in $T$). It remains to show that there is no positive cycle.\\ 
Let $C$ be a positive cycle of $\orient{M}$ and $e$ an edge of $C$. We consider the directed path $P$ of $T$ from the root to $e$ (the root and $e$ included). By definition, the root is in the exterior region of $C$. Let $h$ be the last head of $P$ contained in the exterior region of $C$ and $t$ the tail following $h$ in $P$ (the tail $t$ exists since the last edge $e$ of $P$ is in $C$). By definition, the tail $t$ is either in $C$ or in its interior region. Let $v$ be the end of $h$ (i.e the origin of $t$) in $\orient{M}$ and $h'$  the head of $C$ incident with $v$ (see Figure \ref{fig:dem-merging}). In counterclockwise direction around $v$, we have $h$, $t$ and $h'$ (and possibly some other half-edges). The vertex $v$ is obtained by contracting a cell $C_v$ of the partition $\sim$ corresponding to some vertices of $T$. Each of these vertices is incident to one head in $T$, hence  $h$ and $h'$ were incident to two distinct vertices, say $v_1$ and $v_2$, of $T$. 
The cell $C_v$ is incident to the first corner of $v_1$ which is situated between $h$ and $t$ in counterclockwise direction around $v_1$. Therefore, after the cell $C_v$ is contracted, the half-edges of $v_2$ are situated between $h$ and $t$ in counterclockwise direction around $v$. Thus,  in counterclockwise direction around $v$, we have $h$, $h'$ and $t$ (and possibly some other half-edges). We obtain a contradiction. 
\findem

We now conclude the proof of Theorem \ref{thm:phi-bijective}. \\
\ite Let $\orient{M}$ be a tree-oriented map. We know from Lemma \ref{thm:explosion-produce-tree} that $T=\varphi_0(\orient{M})$ is a tree oriented from the root to the leaves. Moreover, we know from Lemma \ref{thm:explosion-produce-partition} that the partition $\sim=\varphi_0(\orient{M})$ of the vertex set of $T$ is non-crossing. Let $u$ be a vertex of $T$. Let  $\{v_1,\ldots,v_k\}$ be a part of the partition $\sim$ corresponding to a vertex $v$ of $\orient{M}$. The cell $C_v$ associated to $v$ during  the vertex explosion process is incident to the corner of $v_i$, $i=1\ldots k$ at the right of the head incident with $v_i$ (see Figure \ref{fig:vertex-explosion2}). Since  $T$ is oriented from the root to the leaves, this corner is the first corner of $v_i$. Therefore, by definition of $\psi$, we have $\psi\circ\varphi(\orient{M})=\orient{M}$. Thus, $\psi\circ\varphi$ is the identity mapping on  tree-oriented maps.\\
\ite Let $T$ be a tree of size $n$ and $\sim$ be a non-crossing partition on a linearly ordered set $S$ of size $n+1$. We know from Lemma \ref{thm:psi-produces-tree-oriented} that $\orient{M}=\psi(T,\sim)$ is a tree-oriented map. We think to the tree $T$ as being oriented from the root to the leaves and we identify the set $S$ with the vertex set of $T$. Let $v$ be a vertex of $\orient{M}$ corresponding to the part $\{v_1,\ldots,v_k\}$ of the partition $\sim$. The vertex $v$ is obtained by contracting a cell $C_v$ incident with the first corner of $v_i$, $i=1\ldots k$, that is, the corner at the right of the head $h_i$ incident with $v_i$. Therefore, if $t$ is a tail incident with $v_i$ in $T$, then,  $h_i$ is the first head encountered in counterclockwise direction around $v$ starting from $t$ (in $\orient{M}$). Given the definition of the vertex explosion process, the adjacency relations between the half-edges incident with $v$ that are preserved by the vertex explosion process are exactly the adjacency relations in the tree $T$. Thus, the trees $\varphi_0(\orient{M})$ and $T$ are the same. Moreover, the  part of the partition $\varphi_1(\orient{M})$ associated to the vertex $v$ is $\{v_1,\ldots,v_k\}$. Thus, the partitions $\varphi_1(\orient{M})$ and $\sim$ are the same. Hence, $\varphi\circ\psi$ is the identity mapping on pairs made of a tree of size $n$ and a non-crossing partition of size $n+1$.\\
Thus,  the mapping $\varphi$ is a bijection between tree-oriented maps of size $n$ and pairs made of a tree of size $n$ and a non-crossing partition of size $n+1$. This completes the proof of Proposition \ref{thm:varphi-bijective} and Theorem \ref{thm:phi-bijective}. 
\findem

%It is clear that the vertex explosion process and the cell contraction are reverse operations (see Figure \ref{fig:vertex-explosion2}). Therefore, $\varphi$ and $\psi$ are reverse mappings. Thus,  $\varphi$ is a bijection between tree-oriented maps of size $n$ and pairs made of a tree of size $n$ and a non-crossing partition of size $n+1$. This completes the proof of Proposition \ref{thm:varphi-bijective} and Theorem \ref{thm:phi-bijective}. 
%\findem

\section{Correspondence with a bijection due to Cori, Dulucq and Viennot} \label{section:equivalence}
\begin{figure}[ht!]
\begin{center}
\hspace{-.5cm}\input{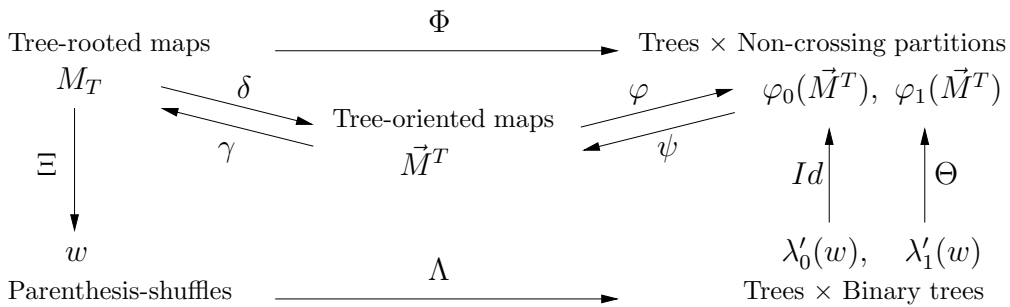}
\caption{The bijection diagram.} \label{fig:equivalence-bijections}
\end{center}
\end{figure}

In this section, we prove that our bijection $\Phi$ is isomorphic to a former bijection due to Cori, Dulucq and Viennot defined on parenthesis-shuffles \cite{Dulucq:shuffle-parenthesis-system}. We know that  tree-rooted maps are in one-to-one correspondence with parenthesis-shuffles by the mapping $\Xi$ defined in Section \ref{section:preliminary}. 
Our bijection $\Phi: M_T \mapsto (\varphi_0(\orientT{M}),\varphi_1(\orientT{M}))$ associates with any tree-rooted map $M_T$ of size $n$, a tree $\varphi_0(\orientT{M})$ of size $n$ and a non-crossing partition $\varphi_1(\orientT{M})$ of size $n+1$. The bijection $\Lambda: w\mapsto(\lambda_0'(w),\lambda_1'(w))$ of Cori \emph{et al.} associates with any parenthesis-shuffle $w$ of size $n$, a tree $\lambda_0'(w)$ of size $n$ and a binary tree $\lambda_1'(w)$ of size $n+1$. 
%Our bijection $\Phi$ associates the ordered pair $(\varphi_0(\orientT{M}),\varphi_1(\orientT{M}))$ with a tree-rooted map $M_T$ of size $n$, where  $\varphi_0(\orientT{M})$ is a tree of size $n$ and $\varphi_1(\orientT{M})$ is a non-crossing partition of size $n+1$. The bijection $\Lambda$ of Cori \emph{et al.} associates the ordered pair $(\lambda_0'(w),\lambda_1'(w))$ with a parenthesis-shuffle $w$ of size $n$, where  $\lambda_0'(w)$ is a tree of size $n$ and $\lambda_1'(w)$ is a binary tree of size $n+1$. 
We shall prove that these two bijections are isomorphic via the encoding of tree-rooted maps by parenthesis-shuffles. That is, we shall prove that there exist two independent bijections $\Omega$ and $\Theta$  such that, if $w=\Xi(M_T)$, then  $\varphi_0(\orientT{M})=\Omega(\lambda_0'(w))$ and $\varphi_1(\orientT{M})=\Theta(\lambda_1'(w))$. In fact, we have adjusted some definitions from \cite{Dulucq:shuffle-parenthesis-system} so that $\Omega$ is the identity mapping on trees. This situation is represented in Figure \ref{fig:equivalence-bijections}.\\

\subsection{The bijection $\Lambda$ of Cori, Dulucq and Viennot}
We begin with a presentation of the bijection $\Lambda$ of Cori \emph{et al.} For the sake of simplicity, the presentation given here is not \emph{completely} identical to the one of the original article \cite{Dulucq:shuffle-parenthesis-system}. But, whenever our definitions differ there is an obvious equivalence via a composition with a simple, well-known bijection. The interested reader can look for more details in the original article. In this article, Cori \emph{et al.} defined recursively two mappings $\lambda_0$ and $\lambda_1$ on the set of \emph{prefix-shuffles}. A prefix-shuffle is a  word $w$ on the alphabet $\{a,\B{a},b,\B{b}\}$ such that, for all prefixes $w'$ of $w$, we have  $|w'|_a\geq |w'|_{\B{a}}$ and  $|w'|_b\geq |w'|_{\B{b}}$. Note that the set of prefix-shuffles is the set of  prefixes of parenthesis-shuffles. The mappings $\lambda_0$ and $\lambda_1$ both eventually return trees. In  the original paper  \cite{Dulucq:shuffle-parenthesis-system}, the trees returned by $\lambda_0$ and $\lambda_1$ were called the \emph{leaf code} and   the \emph{tree code} respectively. \\

\begin{figure}[ht!]
\begin{center}
\input{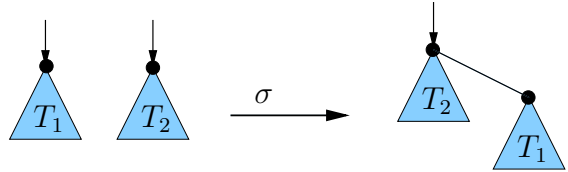}
\caption{The mapping $\sigma$ on ordered pairs of trees.} \label{fig:sigma}
\end{center}
\end{figure}

We first define the mapping $\lambda_0$. It involves the  mapping $\sigma$ that associates the tree $\sigma(T_1,T_2)$ represented in Figure \ref{fig:sigma} with the ordered pair of trees $(T_1,T_2)$.
  
We consider the alphabet $\textbf{U}=\{u,v\}$ and the infinite alphabet $\textbf{T}$ consisting of all trees. A word $s$ on the alphabet  $\textbf{U}\cup \textbf{T}$ is a \emph{tree-sequence} if $s=ut_1u\ldots t_{i-1}u t_i vt_{i+1}\ldots t_k v$ where $1\leq i \leq k$ and $t_1,\ldots,t_k$ are trees. The mapping $\lambda_0$ associates tree-sequences with prefix-shuffles.
\begin{Def}\label{def:lambdaT}
The mapping  $\lambda_0$ is recursively defined on prefix-shuffles by the following rules:
\begin{itemize}
\item If $w=\epsilon$ is the empty word,  $\lambda_0(w)$ is the tree-sequence $u \tau v$ where $\tau$ is the tree reduced to a root and a vertex.\\
 \centerline{\begin{picture}(0,0)%
\includegraphics{rule-T-epsilon.pstex}%
\end{picture}%
\setlength{\unitlength}{789sp}%
\begingroup\makeatletter\ifx\SetFigFont\undefined%
\gdef\SetFigFont#1#2#3#4#5{%
  \reset@font\fontsize{#1}{#2pt}%
  \fontfamily{#3}\fontseries{#4}\fontshape{#5}%
  \selectfont}%
\fi\endgroup%
\begin{picture}(2276,1408)(1501,-2936)
\put(1501,-2311){\makebox(0,0)[lb]{\smash{\SetFigFont{5}{6.0}{\rmdefault}{\mddefault}{\updefault}\large{$\tau~:$}}}}
\end{picture}
}
\item If $w=w'a$, the tree-sequence $\lambda_0(w)$ is obtained from  $\lambda_0(w')$ by replacing the last occurrence of  $u$ by $u\tau v$.
%\centerline{\input{rule-T-a.pstex_t}}
\item If $w=w'b$, the tree-sequence $\lambda_0(w)$ is obtained from  $\lambda_0(w')$ by replacing the first  occurrence of $v$ by $u\tau v$.
%\centerline{\input{rule-T-b.pstex_t}}
\item If $w=w'\B{a}$, we consider the first occurrence of $v$ in $\lambda_0(w')$ and the trees $T_1$ and $T_2$ directly preceding and following it. The tree-sequence $\lambda_0(w)$ is obtained from $\lambda_0(w')$ by  replacing the subword $T_1 v T_2$ by the tree $\sigma(T_1,T_2)$.
%\centerline{\input{rule-T-abar.pstex_t}}
\item If $w=w'\B{b}$,  we consider the last occurrence of $u$ in $\lambda_0(w')$ and the trees $T_1$ and $T_2$ directly preceding and following it. The tree-sequence $\lambda_0(w)$ is obtained from $\lambda_0(w')$ by  replacing  the subword $T_1 u T_2$ by the tree $\sigma(T_1,T_2)$.
%\centerline{\input{rule-T-bbar.pstex_t}}
\end{itemize}
\end{Def}
 
\begin{figure}[ht!]
\begin{center}
\input{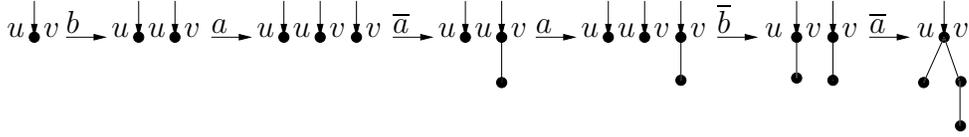}
\vspace{-.7cm} \caption{The mapping $\lambda_0$ applied to the prefix-shuffle $w=ba\B{a}a\B{b}\B{a}$.} \label{fig:exemple-lambda-T}
\end{center}
\end{figure}

We applied the mapping $\lambda_0$ to the word $w=ba\B{a}a\B{b}\B{a}$. The different steps are represented in Figure \ref{fig:exemple-lambda-T}.\\

It is easily seen by induction that the number of $v$ (resp. $u$) in $\lambda_0(w)$ is $|w|_a-|w|_{\B{a}}+1$ (resp. $|w|_b-|w|_{\B{b}}+1$). Hence, the mapping $\lambda_0$ is well defined on prefix-shuffles. Moreover, the first letter $u$ and last letter $v$ are never replaced by anything. Observe also (by induction) that the letters $u$ always precede the letters $v$ in  $\lambda_0(w)$. Thus, $\lambda_0(w)$ is indeed a tree-sequence. If $w$ is a parenthesis-shuffle, there is exactly one letter $u$ and one letter $v$ in $\lambda_0(w)$, hence $\lambda_0(w)$ is a three letter word $uTv$.
\begin{Def}
The mapping $\lambda_0'$ associates with a parenthesis-shuffle $w$ the unique tree $T$ in the tree-sequence $\lambda_0(w)=uTv$.
\end{Def}
Observe that, for any prefix-shuffle $w$, the total number of edges in the trees $t_1,\ldots,t_k$ of the tree-sequence $\lambda_0(w)=ut_1u\ldots t_{i-1}u t_i vt_{i+1}\ldots t_k v$ is $|w|_{\B{a}}+|w|_{\B{b}}$. Hence, if $w$ is parenthesis-shuffle of size $n$, the tree $\lambda_0'(w)$ has size $n$.\\

We now define the mapping $\lambda_1$ which associates \emph{binary trees} with prefix-shuffles. A binary tree is a (planted plane) tree for which each vertex is either of degree 3, a \emph{node}, or of degree 1, a \emph{leaf}. The size of a binary tree is defined as the number of its nodes. It is well-known that binary trees of size $n$ (i.e. with $n$ nodes) are in one-to-one correspondence with trees of size $n$ (i.e. with $n$ edges). \\

\begin{figure}[ht!]
\begin{center}
\input{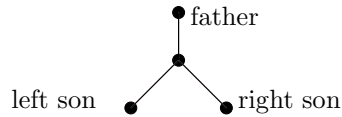}
\caption{Left and right son of a node} \label{fig:left-right-son}
\end{center}
\end{figure}

In a binary tree, the two sons of a node are called \emph{left son} and \emph{right son}. In counterclockwise order around a node we find the father (or the root), the left son and the right son (see Figure \ref{fig:left-right-son}). 
A \emph{left leaf} (resp. \emph{right leaf}) is a leaf which is a left son (resp. right son). As before, we compare vertices according to their order of  appearance around the tree and we shall talk about the \emph{first} and \emph{last} leaf. Moreover, a leaf will be either \emph{active} or \emph{inactive}. Graphically, active leaves will be represented by circles and inactive ones by squares.\\% We are now ready for the definition of the mapping $\lambda_1$.

\begin{Def}\label{def:lambdaP} 
The mapping $\lambda_1$ is recursively defined on prefix-shuffles by the following rules:
\begin{itemize}
\item If $w=\epsilon$ is the empty word,  $\lambda_1(w)$ is the binary tree $B_1$ consisting of a root, a node and two active leaves.\\
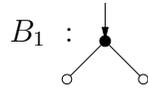
\centerline{\begin{picture}(0,0)%
\includegraphics{rule-P-epsilon.pstex}%
\end{picture}%
\setlength{\unitlength}{789sp}%
\begingroup\makeatletter\ifx\SetFigFont\undefined%
\gdef\SetFigFont#1#2#3#4#5{%
  \reset@font\fontsize{#1}{#2pt}%
  \fontfamily{#3}\fontseries{#4}\fontshape{#5}%
  \selectfont}%
\fi\endgroup%
\begin{picture}(4358,2590)(601,-4118)
\put(601,-2761){\makebox(0,0)[lb]{\smash{{\SetFigFont{5}{6.0}{\rmdefault}{\mddefault}{\updefault}\large{$B_1~:$}}}}}
\end{picture}%
}
\item If $w=w'a$, the tree $\lambda_1(w)$ is obtained from  $\lambda_1(w')$ by replacing the last active left leaf by $B_1$.\\%a node incident to two active leaves.  \\
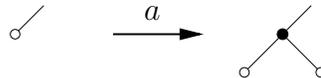
\centerline{\begin{picture}(0,0)%
\includegraphics{rule-P-a.pstex}%
\end{picture}%
\setlength{\unitlength}{789sp}%
\begingroup\makeatletter\ifx\SetFigFont\undefined%
\gdef\SetFigFont#1#2#3#4#5{%
  \reset@font\fontsize{#1}{#2pt}%
  \fontfamily{#3}\fontseries{#4}\fontshape{#5}%
  \selectfont}%
\fi\endgroup%
\begin{picture}(9916,2290)(2243,-4718)
\put(6451,-2911){\makebox(0,0)[lb]{\smash{{\SetFigFont{5}{6.0}{\rmdefault}{\mddefault}{\updefault}\large{$a$}}}}}
\end{picture}%
}
\item If $w=w'b$, the tree $\lambda_1(w)$ is obtained from  $\lambda_1(w')$ by replacing the first active right leaf by $B_1$.\\ %a node incident to two active leaves.  \\
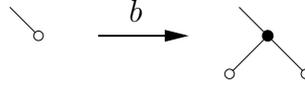
\centerline{\begin{picture}(0,0)%
\includegraphics{rule-P-b.pstex}%
\end{picture}%
\setlength{\unitlength}{789sp}%
\begingroup\makeatletter\ifx\SetFigFont\undefined%
\gdef\SetFigFont#1#2#3#4#5{%
  \reset@font\fontsize{#1}{#2pt}%
  \fontfamily{#3}\fontseries{#4}\fontshape{#5}%
  \selectfont}%
\fi\endgroup%
\begin{picture}(9491,2290)(2668,-4718)
\put(6451,-2911){\makebox(0,0)[lb]{\smash{{\SetFigFont{5}{6.0}{\rmdefault}{\mddefault}{\updefault}\large{$b$}}}}}
\end{picture}%
}
\item If $w=w'\B{a}$, the tree $\lambda_1(w)$ is obtained from  $\lambda_1(w')$ by inactivating the first active right leaf.  \\
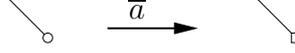
\centerline{\begin{picture}(0,0)%
\includegraphics{rule-P-abar.pstex}%
\end{picture}%
\setlength{\unitlength}{789sp}%
\begingroup\makeatletter\ifx\SetFigFont\undefined%
\gdef\SetFigFont#1#2#3#4#5{%
  \reset@font\fontsize{#1}{#2pt}%
  \fontfamily{#3}\fontseries{#4}\fontshape{#5}%
  \selectfont}%
\fi\endgroup%
\begin{picture}(9195,1395)(2368,-4123)
\put(6151,-3361){\makebox(0,0)[lb]{\smash{{\SetFigFont{5}{6.0}{\rmdefault}{\mddefault}{\updefault}\large{$\B{a}$}}}}}
\end{picture}%
}
\item If $w=w'\B{b}$, the tree $\lambda_1(w)$ is obtained from  $\lambda_1(w')$ by inactivating the  last active left leaf. \\ 
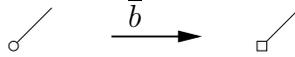
\centerline{\begin{picture}(0,0)%
\includegraphics{rule-P-bbar.pstex}%
\end{picture}%
\setlength{\unitlength}{789sp}%
\begingroup\makeatletter\ifx\SetFigFont\undefined%
\gdef\SetFigFont#1#2#3#4#5{%
  \reset@font\fontsize{#1}{#2pt}%
  \fontfamily{#3}\fontseries{#4}\fontshape{#5}%
  \selectfont}%
\fi\endgroup%
\begin{picture}(9191,1395)(2243,-4123)
\put(6001,-3361){\makebox(0,0)[lb]{\smash{{\SetFigFont{5}{6.0}{\rmdefault}{\mddefault}{\updefault}\large{$\B{b}$}}}}}
\end{picture}%
}
\end{itemize}
\end{Def}

\begin{figure}[ht!]
\begin{center}
\input{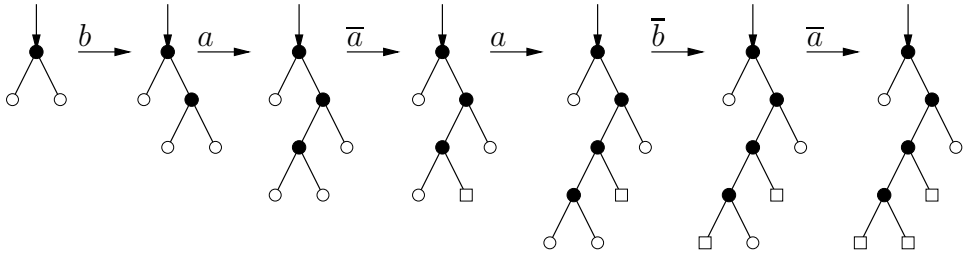}
\caption{The mapping $\lambda_1$ on the word $w=ba\B{a}a\B{b}\B{a}$.} \label{fig:exemple-lambda-P}
\end{center}
\end{figure}

 We applied the mapping $\lambda_1$ to the word $w=ba\B{a}a\B{b}\B{a}$. The different steps are represented in Figure \ref{fig:exemple-lambda-P}.

It is easily seen by induction that the number of active right leaves (resp. left leaves) in $\lambda_1(w)$ is $|w|_a-|w|_{\B{a}}+1$ (resp. $|w|_b-|w|_{\B{b}}+1$). Hence, the mapping $\lambda_1$ is well defined on prefix-shuffles. Observe that the binary tree $\lambda_1(w)$ has $|w|_a+|w|_b+1$ nodes. Observe  also (by induction) that active left leaves always precede active right leaves in  $\lambda_1(w)$. 
Moreover, if $w$ is a parenthesis-shuffle, only the first left leaf and the last right leaf are active (since they can never be inactivated). 
\begin{Def}
The mapping $\lambda_1'$ associates with a parenthesis-shuffle $w$ of size $n$ the binary tree of size $n+1$ obtained from $\lambda_1(w)$ by inactivating the two active leaves.
\end{Def}

We now make some informal remarks explaining why the mapping $w \mapsto (\lambda_0(w),\lambda_1(w))$ is injective. It is, of course, possible to decide from $(\lambda_0(w),\lambda_1(w))$ if  $w$ is the empty word. Indeed, $w$ is the empty word iff $\lambda_1(w)=B_1$ (equivalently iff $\lambda_0(w)=\tau$). Otherwise, the remarks below show that the last letter $\alpha$ of $w=w'\alpha$ can be determined as well as   $\lambda_0(w')$ and $\lambda_1(w')$. So any prefix-shuffle $w$ can be entirely recovered from $(\lambda_0(w),\lambda_1(w))$.\\

\noindent \textbf{Remarks:}\\
\ite For any prefix-shuffle $w$, the number of letters $u$ (resp. $v$) in  the tree-sequence $\lambda_0(w)$ is equal to the number of active left leaves (resp. right leaves) in the binary tree $\lambda_1(w)$.  Furthermore, it can be shown by induction that the size of the tree $t_i$ lying between the $i^{th}$ and $i+1^{th}$ letters  $u,v$ in $\lambda_0(w)$ is the number of inactive leaves lying between the $i^{th}$ and $i+1^{th}$ active leaves in $\lambda_1(w)$. 

\ite The three following statements are equivalent:\\
\iten the word $w$ is not empty and the last letter $\alpha$ of $w=w'\alpha$ is in $\{a,b\}$, \\
\iten there is a sequence $u\tau v$ in $\lambda_0(w)$,\\
\iten there is an active left leaf and an active right leaf which are siblings. \\
In this case, $\lambda_1(w')$ is obtained from $\lambda_1(w)$ by deleting the two actives leaves and making the father an active leaf $\ell$. Moreover, $\alpha=a$ (resp. $\alpha=b$) if  $\ell$ is a left leaf (resp. right leaf) in $\lambda_1(w')$ in which case  $\lambda_0(w')$ is obtained from $\lambda_0(w)$ by replacing the subword $u\tau v$ by $u$ (resp. $v$). 

\ite If the last letter $\alpha$ of $w=w'\alpha$ is in $\{\B{a},\B{b}\}$, we know from the above remark that the tree $T$ lying between the last letter $u$ and the first letter $v$ in the tree-sequence $\lambda_0(w)$ has size $k>0$. Since $k>0$, the tree  $T$ admits a (unique) preimage $(T_1,T_2)$ by the mapping $\sigma$. Let $k'$ be the size of the tree $T_1$. Then $k'<k$. We know that there are $k$ inactive leaves lying between the last active left leaf and the first active right leaf in $\lambda_1(w)$.  The binary tree $\lambda_1(w')$ is obtained from $\lambda_1(w)$ by activating the  $k'+1^{th}$ leaf $\ell$ encountered  when following the border of the tree starting from the last active left leaf. Moreover, $\alpha=\B{a}$ (resp. $\alpha=\B{b}$) if $\ell$ is a right leaf (resp. left leaf), in which case the tree-sequence $\lambda_0(w')$ is obtained from $\lambda_0(w)$ by replacing $T$ by $T_1 v T_2$ (resp. $T_1 u T_2$).\\

From these remarks, we see that the mapping $w \mapsto (\lambda_0(w),\lambda_1(w))$ is injective. It can be shown, with the same ideas, that it is bijective on the set of pairs consisting of a tree-sequence $S$ and a binary tree $B$ with active and inactive leaves satisfying the  following conditions:\\ 
%\begin{enumerate}
\iten the active left leaves precede the active right leaves in $B$,\\
\iten the number of active left leaves (resp. right leaves) in $B$ is the same as the number of $u$ (resp. $v$) in $S$,\\
\iten the number of inactive leaves lying between the $i^{th}$ and $i+1^{th}$ active leaves in $B$ is the size of the tree lying between the $i^{th}$ and $i+1^{th}$ letters $u,v$ in $S$.\\
%\end{enumerate}

We now define the mapping $\Lambda$ of Cori \emph{et al.} on parenthesis-shuffles.
\begin{Def} 
The mapping $w\mapsto (\lambda_0'(w),\lambda_1'(w))$ defined on parenthesis-shuffles is denoted $\Lambda$.
\end{Def}
We know that $\Lambda$ associates with a parenthesis-shuffle of size $n$ a pair consisting of a tree of size $n$ and a binary tree of size $n+1$. The remarks above should convince the reader that the mapping $\Lambda$ is a bijection between these two sets of objects.\\

\subsection{The bijections $\Phi$ and $\Lambda$ are isomorphic}
We now return to our business and prove that  the bijection $\Lambda$ of Cori \emph{et al.} and  our bijection $\Phi$ 
%presented in Section \ref{section:bijection} 
are isomorphic. Before stating precisely this result, we define a (non-classical) bijection $\theta$ between binary trees and trees. By composition, this allows us to define a  bijection $\Theta$ between binary trees and non-crossing partitions.\\

\begin{figure}[ht!]
\begin{center}
\input{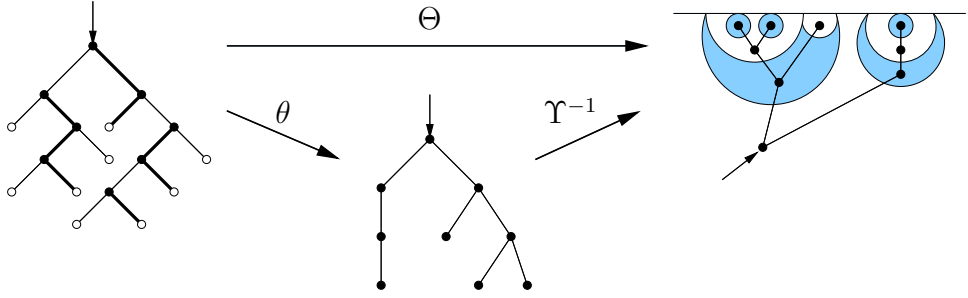}
\caption{The mappings $\theta$ and $\Theta$.} \label{fig:exemple-theta}
\end{center}
\end{figure}

Let $e$ be an edge of a binary tree. The edge $e$ is said to be branching if one of its vertices is a right son and the other is a left son or the root-vertex. Intuitively, this means that the edge $e$ is non-parallel to its parent-edge. For instance, the branching edges of the binary tree in Figure \ref{fig:exemple-theta} are indicated by thick lines.

\begin{Def} \label{def:theta}
Let $B$ be a binary tree. The tree $\theta(B)$ is obtained by contracting every non-branching edge. The non-crossing partition  $\Theta(B)$ is the image of  $\theta(B)$ by the mapping $\Upsilon^{-1}$ (see  Figure \ref{fig:partition-tree}).
\end{Def}

We applied the mapping $\Theta$ to the binary tree of Figure  \ref{fig:exemple-theta}.

%Consider a binary tree $B$. Changing a leaf  into a node with two leaves increases by one the number of branching edges. It follows by induction that the size of the tree $\theta(B)$ (i.e. its number of edges) is the size of the binary tree $B$ (i.e. its number of nodes). Therefore, the size of the non-crossing partition $\Theta(B)$ is the size of $B$.\\
%We now state the main result of this section.

The mapping $\Theta$ is a bijection between binary trees of size $n$ ($n$ nodes) and trees of size $n$ ($n$ edges). The proof is omitted here since we will not use this property. \\
We now state the main result of this section.

\begin{thm}\label{thm:equivalence}
Let $M_T$ be a tree-rooted map and $w=\Xi(M_T)$ its associated parenthesis-shuffle. Let $\varphi_0(\orientT{M})$ and   $\varphi_1(\orientT{M})$ be the tree and the non-crossing partition obtained from $M_T$ by the mapping $\Phi$. Let $\lambda_0'(w)$ and $\lambda_1'(w)$  be the tree and binary tree obtained from $w$ by the mapping $\Lambda$.
Then $\varphi_0(\orientT{M})=\lambda_0'(w)$ and   $\varphi_1(\orientT{M})=\Theta(\lambda_1'(w))$. 
\end{thm}

\begin{figure}[ht!]
\begin{center}
\input{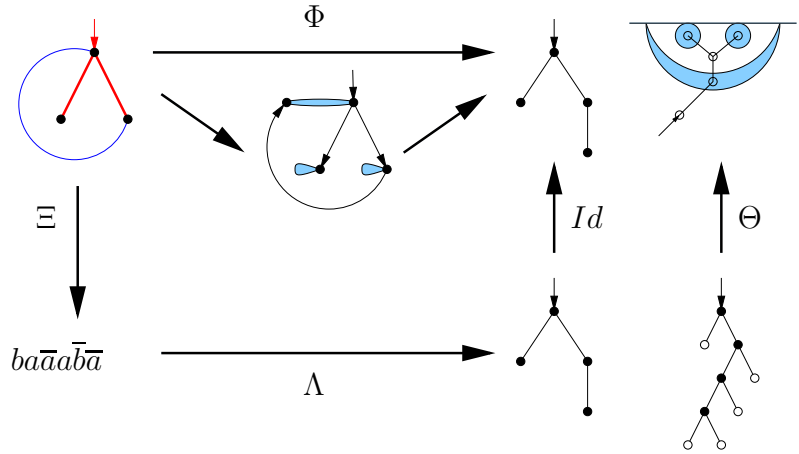}
\caption{The isomorphism between $\Lambda$ and $\Phi$.} \label{fig:exp-equivalence}
\end{center}
\end{figure}

This relation between the mappings  $\Lambda$ and $\Phi$ is represented by Figure \ref{fig:equivalence-bijections}. As an illustration, we applied the mapping $\Phi$ to the tree-rooted map $M_T$ of Figure \ref{fig:exp-equivalence} and we applied the mapping $\Lambda$ to $w=\Xi(M_T)=ba\B{a}a\B{b}\B{a}$. The rest of this section is devoted to the proof of Theorem \ref{thm:equivalence}.\\

\subsection{Prefix-maps} \label{section:prefix-maps}
The mappings $\lambda_0'$ and $\lambda_1'$ are defined on parenthesis-shuffles from the more general mappings  $\lambda_0$ and $\lambda_1$ defined on prefix-shuffles. In order to relate $\varphi_0(\orientT{M})$ and $\lambda_0'(w)$ (resp. $\varphi_1(\orientT{M})$ and $\lambda_1'(w)$)  we need to define the \emph{prefix-maps} which are in one-to-one correspondence with prefix-shuffles. As we will see, prefix-maps are tree-oriented maps together with some \emph{dangling} heads in the root-face. 
%Moreover, in prefix-maps, some edges are said to be \emph{active}. 
In Subsections \ref{section:varphi0=lambda0} and \ref{section:varphi1=lambda1}  we shall extend the mappings $\varphi_0$ and $\varphi_1$ defined in Section \ref{section:bijection} to prefix-maps.\\

For any prefix-shuffle $w$ we denote by $w_a$ (resp. $w_b$) the subword of $w$ consisting of the letters $a,\B{a}$ (resp. $b,\B{b}$). The words $w_a$ and $w_b$ are prefixes of parenthesis systems. We say that an occurrence of a letter $c=a,b$ is \emph{paired} with an occurrence of $\B{c}$ if the subword of $w_c$ lying between these two letters is a parenthesis system. There are $|w|_a -|w|_{\B{a}}$ non-paired letters $a$  and  $|w|_b -|w|_{\B{b}}$ non-paired letters $b$ in $w$.  We denote by $w_a^{+}$  the parenthesis system obtained from $w_a$ by adding $|w|_a -|w|_{\B{a}}$ letters $\B{a}$ at the end of this word. \\
%We denote by $w_a^{+}$  (resp. $w^+$) the word $w_a$ (resp. $w$) followed by  $|w|_a -|w|_{\B{a}}$ letters $\B{a}$.  Note that $w_a^{+}$ is a parenthesis system.

\begin{figure}[h!]
\begin{center}
\input{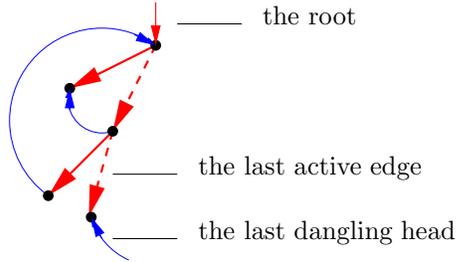}
\caption{The prefix-map associated to $bab\B{a}a\B{b}a\B{b}\B{a}ab$.} \label{fig:exp-carte-prefixe2}
\end{center}
\end{figure}

Let $w$ be a prefix-shuffle. We define  $T_w$ as the tree associated to the parenthesis system $w_a^{+}$, that is, $T_w$ is such that, making the tour of $T_w$  and writing  $a$ the first time we follow an edge and $\B{a}$  the second time, we obtain  $w_a^{+}$. We orient the edges of $T_w$ from the root to the leaves. Then, we add half-edges to $T_w$ by looking at the position of the letters $b$ and $\B{b}$ in $w$. More precisely, we read the word $w$ and while making the tour of $T$ according to the letters $a,\B{a}$, we insert heads for the letters $b$ and tails for the letters $\B{b}$. If an occurrence of $b$ and an occurrence of $\B{b}$ are paired in $w$ we connect the corresponding head and tail. We obtain an oriented map together with some heads called \emph{dangling heads} corresponding to non-paired letters $b$ of $w$. In the tree $T_w$, the edges corresponding to non-paired letters $a$ are called \emph{active} while the others are called \emph{inactive}. We denote by $M_w$, and call \emph{prefix-map associated with $w$}, the oriented map (with dangling heads and active edges) obtained. For instance, the prefix-map associated with $bab\B{a}a\B{b}a\B{b}\B{a}ab$ has been represented in Figure \ref{fig:exp-carte-prefixe2} (the active edges are dashed). \\ 
Observe that $T_w$ is a spanning tree of the prefix-map $M_w$. %In the following, $T_w$ is called  the \emph{spanning tree of the prefix-map $M_w$}. 
The orientation of $M_w$ is the tree-orientation associated to the spanning tree $T_w$ by the mapping $\delta$ defined in Section \ref{section:bijection}. In particular, when $w$ is a parenthesis-shuffle, the prefix-map $M_w$ is a map (i.e. it has no active edge and no dangling head except for the root) which is tree-oriented. More precisely, if $w=\Xi(M_T)$, the tree-oriented map $M_w$ is $\orientT{M}\equiv\delta(M_T)$.\\

%Let $w$ be a prefix-shuffle. We define  $T_w$ as the tree associated to the parenthesis system $w_a^{+}$, that is, $T_w$ is such that, making the tour of $T_w$  and writing  $a$ the first time we follow an edge and $\B{a}$  the second time, we obtain  $w_a^{+}$. We orient the edges of $T_w$ from the root to the leaves. Then, we add half-edges to $T_w$ by looking at the position of the letters $b$ and $\B{b}$ in $w$. More precisely, we read the word $w$ and while making the tour of $T$ according to the letters $a,\B{a}$, we insert heads for the letters $b$ and tails for the letters $\B{b}$. If an occurrence of $b$ and an occurrence of $\B{b}$ are paired in $w$ we connect the corresponding head and tail.We denote by $M_w$, and call \emph{prefix-map associated with $w$}, the map  obtained. The tree $T_w$ is called the \emph{spanning tree of the prefix-map $M_w$}. For instance, the prefix-map associated with $bab\B{a}a\B{b}a\B{b}\B{a}ab$ has been represented in Figure \ref{fig:exp-carte-prefixe2} (forget that some edges are dashed for the time being). \\
%Observe that, if $w$ is a parenthesis-shuffle, the prefix-map $M_w$ is a map (i.e. it has no active edge and no dangling-head except for the root) an it is tree-oriented. More precisely, if $w=\Xi(M_T)$, the tree-oriented map $M_w$ is $\orientT{M}$.\\

%In the spanning tree $T_w$, edges corresponding to non-paired letters $a$ are called \emph{active} while the others are called \emph{inactive}. For instance, the active edges are the two dashed edges in Figure \ref{fig:exp-carte-prefixe2}. 
Let $w$ be a prefix-shuffle. The heads of active edges in the prefix map $M_w$ are called  \emph{rooting heads}, and their ends are called  \emph{rooting vertices}. By convention, the root is considered as a rooting head. As before, we compare active edges (resp. rooting vertices, dangling heads) of $M_w$ according to their order of appearance around $T_w$. By convention, the root is considered as the first rooting head. \\
%In prefix-maps, the heads distinct from the root that are not part of an edge  correspond to non-paired letters $b$ of $w$. These heads are said \emph{dangling}.  We compare active edges (resp. rooting vertices, dangling heads) according to their order of appearance around $T_w$. By convention, the root is considered as the first rooting head. \\
Let $w^{+}$ be the word $w$ followed by $|w|_a -|w|_{\B{a}}$ letters $\B{a}$. We obtain $w^{+}$ by making the tour of the tree $T_w$ and writing $a$ the first time we follow an edge of the tree, $\B{a}$ the second time, $b$ when we cross a head not in the tree and $\B{b}$ when we cross a tail not in the tree. Each prefix of $w^+$ corresponds to a given time in this journey. In particular, $w$ corresponds to a given corner $c$ of a vertex $v$. The $|w|_a -|w|_{\B{a}}$ letters $\B{a}$ at the end of  $w^{+}$ correspond to the left border of active edges followed from $c$ to the root. Thus, the active edges are the edges on the directed path of $T_w$ from the root to $v$. Note that an active edge precedes another one if it appears before on the path from the root to $v$. Therefore, $v$ is the last rooting vertex and $c$ is the corner at the left of the last rooting head. Moreover, active edges are directed from a rooting vertex to the next one (for the appearance order). In particular, the next-to-last rooting vertex (if it exists) is the origin of the last active edge.\\

We now  explore the relation between $M_w$ and $M_{w\alpha}$ when $\alpha$ is a letter in $\{a,\B{a},b,\B{b}\}$. 
\begin{lemma} \label{thm:prefix-map-evolution}
Let $c$ be the corner at the left of the last rooting head of $M_w$. 
\begin{itemize}
\item $M_{wa}$ is obtained from $M_w$ by adding an edge $e$ in the corner $c$. It is oriented from this corner to a vertex not present in $M_w$. The edge $e$ is the  last active edge of $M_{wa}$.
\item $M_{wb}$ is obtained from $M_w$ by adding a dangling head $h$ in the corner $c$. The head $h$ is the last dangling head of $M_{wb}$.
\item $M_{w\B{a}}$ is obtained from $M_w$ by inactivating the last active edge $e$. The origin of $e$ becomes the last rooting vertex.
\item $M_{w\B{b}}$ is obtained from $M_w$ by adding a tail in the corner $c$ and connecting it to the last dangling head.
\end{itemize}
In any case, the appearance order on the edges, half-edges and vertices present in $M_w$ is the same in $M_{w\alpha}$.
\end{lemma}

\dem
As mentioned above, the corner $c$ is the corner reached when the word $w$ is written during the tour of $T_w$ in $M_w$.\\
%\begin{itemize}
\ite Case $\alpha=a$. The letter $a$ added to $w$ is not paired. Therefore, it corresponds to a new active edge $e$ added to $T_w$. This new edge is added in the corner $c$.  The edge $e$ is oriented from $c$ to a new vertex (since it is leaf of $T_{wa}$). All active edges of $M_w$ are encountered before $c$ around the spanning tree $T_w$.  Therefore, $e$ is the last active edge of $M_{wa}$.\\
\ite Case $\alpha=b$. The letter $b$ added to $w$ is not paired.  Therefore, it corresponds to a new dangling head $h$.  This new head is added in the corner $c$.  All dangling heads of $M_w$ are encountered before $c$ around the spanning tree $T_w$. Therefore, $h$  is the last dangling head of $M_{wb}$.\\
\ite Case $\alpha=\B{a}$. The last letter $a$ of $w$ is paired with the letter $\B{a}$ added to $w$. This last letter $a$ corresponds to the last active edge. Therefore, the last active edge $e$ of $M_w$ is inactivated. We know that the next-to-last rooting vertex of $M_w$ is the origin $v$ of the last active edge $e$. Therefore, $v$ becomes the last rooting vertex.\\ 
\ite Case $\alpha=\B{b}$. The last letter $b$ of $w$ is paired with the letter $\B{b}$ added to $w$. This last letter $b$ corresponds to the last dangling head $h'$. Hence, $M_{w\B{b}}$ is obtained from $M_w$ by adding a tail $h$ in the corner $c$ and connecting it to $h'$.
%\end{itemize}
\findem

\noindent This completes our study of prefix-maps. We are now ready to extend the mappings $\varphi_0$ and $\varphi_1$ to prefix maps and to prove  Theorem  \ref{thm:equivalence}.

\subsection{The trees $\varphi_0(\vec{M}^T)$ and $\lambda_0'(w)$ are the same} \label{section:varphi0=lambda0}
In this subsection, we prove that, when $w=\Xi(M_T)$, the trees $\varphi_0(\orientT{M})$ and $\lambda_0'(w)$ are the same. \\

\begin{figure}[h!]
\begin{center}
\input{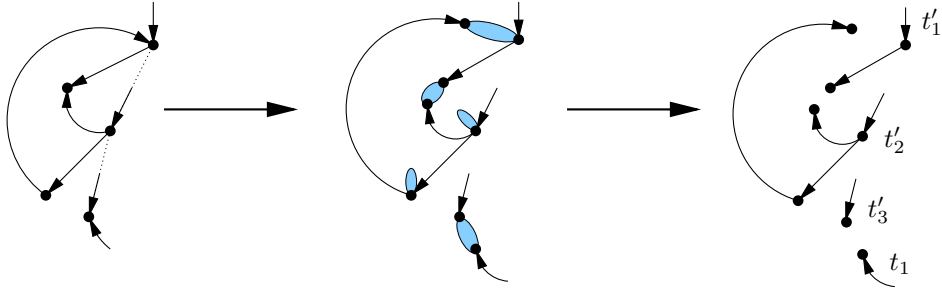}
\caption{The prefix-forest $F_w$ (on the right).} \label{fig:explose-carte-prefixe}
\end{center}
\end{figure}

Let $w$ be a prefix-shuffle and $M_w$ the corresponding prefix-map. Note that any vertex of $M_w$ is incident to at least one head. The \emph{prefix-forest} of $w$, denoted by $F_w$, is obtained by deleting the tails of active edges and then applying the vertex explosion process of Figure \ref{fig:vertex-explosion2} (we forget about the cells corresponding to the parts of the non-crossing partition). We will prove that the prefix-forest is indeed a forest (i.e. a collection of trees) in Proposition \ref{thm:equiv-lambda0-varphi0}. For instance, we represented the prefix-forest of $w=bab\B{a}a\B{b}a\B{b}\B{a}ab$ in Figure \ref{fig:explose-carte-prefixe}. \\

Note that, if $w=\Xi(M_T)$ is a parenthesis-shuffle, the prefix-map $M_w$ is $\orientT{M}$ and no edge is active. Thus, in this case, the prefix-forest $F_w$ is the tree $\varphi_0(\orientT{M})$. We now prove a relation between the prefix-forest $F_w$ and the tree-sequence $\lambda_0(w)$.
\begin{prop} \label{thm:equiv-lambda0-varphi0}
Let $w$ be a prefix-shuffle. Let $h_1<\cdots<h_k$ be the dangling heads and $h_1'<\cdots<h_l'$ be the rooting heads of the prefix-map $M_w$ (linearly ordered by the appearance order).  The prefix-forest $F_w$ is a collection of $k+l$ trees $t_1,\ldots,t_k,t_1',\ldots,t_l'$.  The root of the tree $t_i,i=1,\ldots,k$ is $h_i$ and the root of the tree $t_i',i=1,\ldots,l$ is $h_i'$. Moreover, the tree-sequence $\lambda_0(w)$ is  $u t_1 u\ldots u t_k u t_l' v \ldots v t_1' v~.$
\end{prop}

\begin{figure}[ht!]
\begin{center}
\input{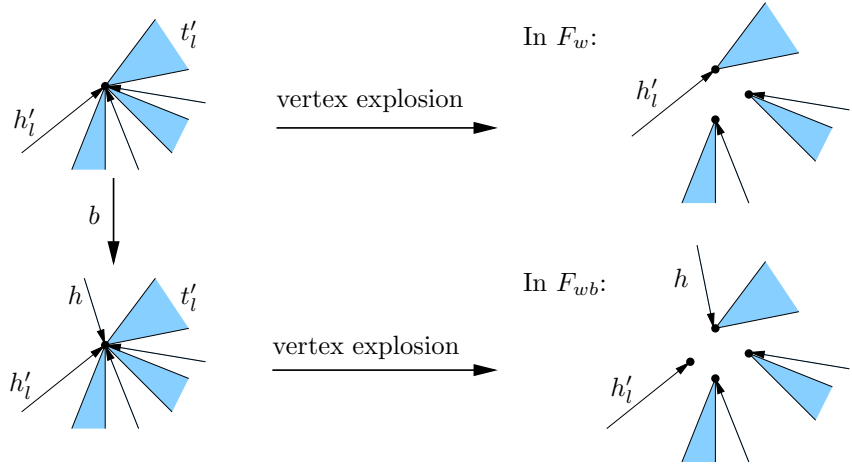}
\caption{The case $\alpha=b$.} \label{fig:dem-equivT}
\end{center}
\end{figure}

\dem
We use Lemma \ref{thm:prefix-map-evolution} and prove the property by induction on the length of $w$.\\ 
If $w$ is the empty word, the prefix-map $M_w$ is the tree $\tau$ reduced to a vertex and a root. Hence, the prefix-forest $F_w$ is reduced to a single tree $\tau=t_1'$. The tree-sequence $\lambda_0(w')$ is equal to $u\tau v$ thus the property is satisfied. If $w'=w\alpha$, we suppose the lemma true for $w$, we write $\lambda_0(w)=u t_1 u\ldots u t_k u t_l' v \ldots v t_1' v$ and study separately the four possible cases.\\
%\begin{itemize}
\ite Case $\alpha=a$. The prefix-map $M_{wa}$ is obtained from $M_w$ by adding an edge $e$ incident to the last rooting vertex. The edge $e$ is the  last active edge of $M_{wa}$. It is oriented toward a new vertex $v$ not present in $M_w$. The tail of $e$ is deleted in the construction of $F_{wa}$ and its head $h=h_{l+1}'$  is only incident to $v$. Therefore,  $F_{wa}$ is obtained from $F_w$ by adding the tree  $\tau=t_{l+1}'$ (the tree reduced to a root and a vertex) rooted on the last rooting head $h$.\\
By definition,  $\lambda_0(wa)=u t_1 u\ldots u t_k u \tau v t_l'v \ldots v t_1' v$, so we observe that the property is satisfied by $wa$.\\
\ite Case $\alpha=b$. The prefix-map $M_{wb}$ is obtained from $M_w$ by adding a dangling head $h=h_{k+1}$ in the corner at the left of the last rooting head $h_l'$. Therefore, during the vertex explosion process $h$ ''steals'' the tree $t_l'$ rooted on $h_l'$ in $F_w$ (see Figure \ref{fig:dem-equivT}). That is, in $F_{wb}$ the tree rooted on $h_l'$ is reduced to a vertex and the tree rooted on  $h$ is $t_l'$.  The head $h$ is the last dangling head of $M_{wb}$. 

\begin{figure}[h!]
\begin{center}
\input{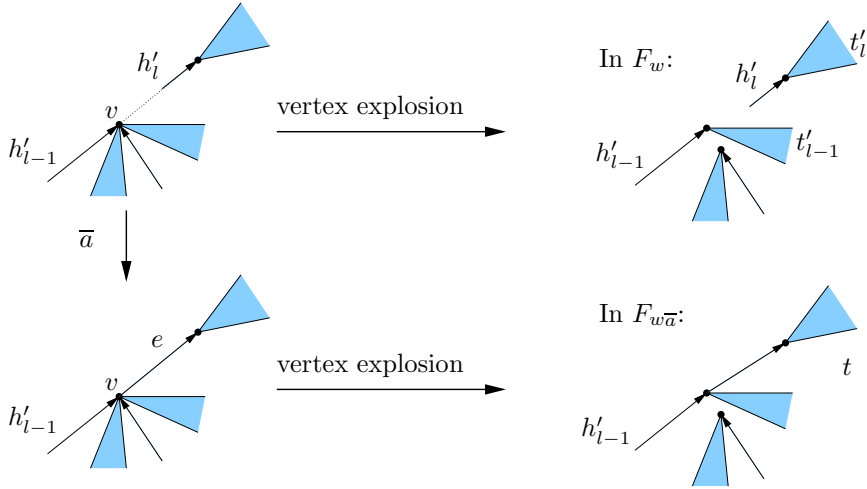}
\caption{The case $\alpha=\B{a}$.} \label{fig:dem-equivT2}
\end{center}
\end{figure}

\noindent  By definition,  $\lambda_0(wb)=u t_1 u\ldots u t_k u t_l' u \tau v t_{l-1}' \ldots v t_1' v,$ so we observe that the property is satisfied by $wb$.\\
\ite Case $\alpha=\B{a}$. The prefix-map  $M_{w\B{a}}$ is obtained from $M_w$ by inactivating the last active edge $e$. The origin $v$ of $e$ is the next-to-last rooting vertex of $M_w$. Moreover, $e$ is the first edge encountered in clockwise order around $v$ starting from $h_{l-1}'$. In $F_{w\B{a}}$, the head $h_l'$ is part of the edge $e$ which links the tree $t_l'$ to the tree $t_{l-1}'$ rooted on $h_{l-1}'$ (see Figure \ref{fig:dem-equivT2}). Therefore, the tree rooted on $h_{l-1}'$ in $F_{w\B{a}}$ is $t=\sigma(t_l',t_{l-1}')$.

\noindent  By definition,  $\lambda_0(w\B{a})=u t_1 u\ldots u t_k u t v t_{l-2}' \ldots v t_1' v,$ so we observe that the property is satisfied by $w\B{a}$.\\
\ite Case $\alpha=\B{b}$.  The prefix-map $M_{w\B{b}}$ is obtained from $M_w$ by adding a tail in the corner at the left of the last rooting head $h_l'$ and connecting it to the last dangling head $h_k$.
In $F_{w\B{b}}$, the head $h_k$ is part of an edge $e$ which links the tree $t_k$ to the tree $t_{l}'$ rooted on $h_{l}'$. Therefore, the tree rooted on $h_{l}'$ in $F_{w\B{b}}$ is $t=\sigma(t_k,t_{l}')$.  The illustration would be the same as Figure \ref{fig:dem-equivT2} except $h_{l-1}',h_l',t_{l-1}',t_l'$ would be replaced by $h_{l}',h_k,t_{l}',t_k$ respectively.\\
\noindent By definition,   $\lambda_0(w\B{b})=u t_1 u\ldots t_{k-1} u t v t_{l-1}' v \ldots v t_1' v,$ so we observe that the property is satisfied by $w\B{b}$.
%\end{itemize}
\findem

As mentioned above, when $w$ is a parenthesis-shuffle $w=\Xi(M_T)$, the prefix-map $M_w$ is the tree-oriented map $\orientT{M}$ and the prefix-forest $F_w$ is the tree $\varphi_0(\orientT{M})$. Therefore, Proposition \ref{thm:equiv-lambda0-varphi0} implies that the tree-sequence $\lambda_0(w)$ is equal to  $u \varphi_0(\orientT{M}) v$. Thus, the trees $\lambda_0'(w)$ and $\varphi_0(\orientT{M})$ are the same.
\findem

\subsection{The partitions $\varphi_1(\vec{M}^T)$ and $\Theta\circ \lambda_1'(w)$ are the same}\label{section:varphi1=lambda1}
In this subsection, we prove that, when $w=\Xi(M_T)$, the non-crossing partition $\varphi_1(\orientT{M})$ is the image of the binary tree $\lambda_1'(w)$ by the mapping $\Theta$ defined in Definition \ref{def:theta}. \\

\begin{figure}[ht!]
\begin{center}
\input{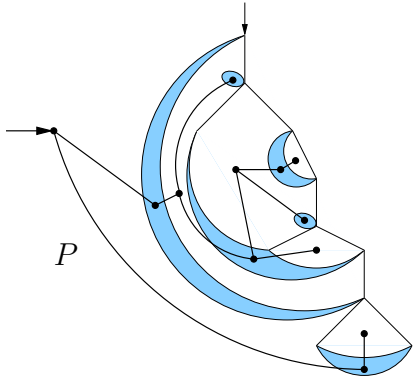}
\caption{The partition-tree of a tree-rooted map.} \label{fig:direct-partition-tree}
\end{center}
\end{figure}

Let $M_T$ be a tree-rooted map. We call \emph{partition-tree} of $M_T$ the tree  $P=\Upsilon\circ \varphi_1(\orientT{M})$ (the mapping $\Upsilon$ is represented in Figure \ref{fig:partition-tree}).  Observe that the tree $P$ can be drawn directly  on the map obtained after the vertex explosion process of Figure \ref{fig:vertex-explosion2}. To do so, one keeps the cells corresponding to the vertices of $\orientT{M}$/ (These cells are glued to the first corner of the vertices of the tree  $\varphi_0(\orientT{M})$). Then, one draws a vertex in each face of $M_T$ and in each cell corresponding to a vertex of $M_T$: this gives the vertices of $P$. The edges of $P$ join vertices in adjacent cells and faces. The tree is rooted canonically. In particular, the root-vertex of $P$ lies in the root-face of $M_T$. This construction is illustrated in Figure \ref{fig:direct-partition-tree}.\\

\begin{figure}[ht!]
\begin{center}
\input{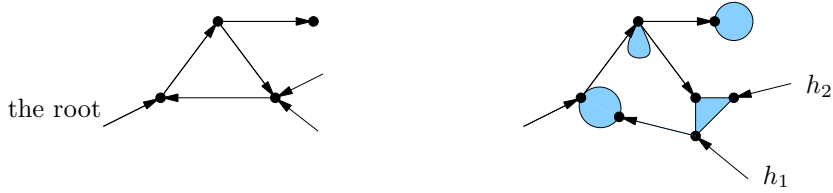}
\caption{The prefix-map associated to $w=baa\B{b}bb\B{a}a$ and the map $\Ms_w$. } \label{fig:partition-tree1}
\end{center}
\end{figure}

We want to extend this construction to prefix-maps. We need some extra vocabulary. Consider a prefix-shuffle $w$ and the corresponding prefix map $M_w$. 
We denote by $\Ms_w$ the map obtained after the vertex explosion process when one keeps the cells corresponding to the vertices of $M_w$.  A face of $\Ms_w$  is said \emph{white} if it corresponds to a face of $M_w$ and \emph{black} if it corresponds to a vertex of $M_w$. For instance, the map $\Ms_w$ in Figure \ref{fig:partition-tree1} has 2 white faces and 4 black faces. We call \emph{regular} the edges of  $M_w$, and \emph{permeable} the edges that separate black and white faces.   The map $\Ms_w$ inherits the root of $M_w$. In particular, it has the same root-face. 
%The permeable edges which are said \emph{external} if they are incident to the root-face and \emph{internal} otherwise.   
The map $\Ms_w$ has $k=|w|_b-|w|_{\B{b}}$ dangling heads which are all in the root-face. We can compare these heads according to their order of appearance \emph{around} the root-face, that is, when following its border in counterclockwise direction starting from the root. We denote by $h_1,\ldots,h_k$ the heads of $\Ms_w$ encountered in this order around the root-face.
%We denote by $h_1,\ldots,h_k$ the $k=|w|_b-|w|_{\B{b}}$ dangling heads of $M_w$. Let us make the \emph{tour of the root-face}, that is, we follow its border in counterclockwise direction starting from the root. We denote by $h_1,\ldots,h_k$ the $k=|w|_b-|w|_{\B{b}}$ dangling heads of $M_w$ encountered in that order during the tour of the root-face.

%We denote by $h_1,\ldots,h_k$ the $k=|w|_b-|w|_{\B{b}}$ dangling heads of $M_w$. Note that the dangling edges are all in the root-face.

\begin{figure}[ht!]
\begin{center}
\input{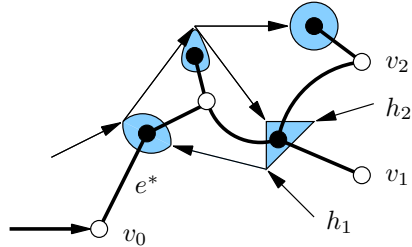}
\caption{The partition-tree $P_w$ (thick lines) drawn on $\Ms_w$ ($w=baa\B{b}bb\B{a}a$).} \label{fig:partition-tree2}
\end{center}
\end{figure}

We define the \emph{partition-tree} $P_w$ of the prefix-map $M_w$ as follows. (We shall prove later that the partition-tree is indeed a tree.) We draw a vertex in each face of $\Ms_w$. The vertex $v_0$ drawn in the root-face is called the \emph{exterior vertex}. We draw $k$ additional vertices  $v_1,\ldots,v_k$ in the root-face, each associated to a dangling head ($v_i$ is associated to $h_i$). These are the vertices of $P_w$. The edges of $P_w$ are the \emph{duals} of permeable edges. 
We need to be more precise. If $e$ is a permeable edge that is not incident to the root-face, its dual joins the vertices drawn in the incident black and white faces. If $e$ is a permeable edge incident to the root-face and a black face $f$, its dual joins the vertex drawn in $f$ to $v_i$ if $h_i$ is the last dangling head encountered before $e$ around the root-face, or to $v_0$ if no dangling head precedes $e$.
%For internal permeable edges this dual is well defined. (If  an internal permeable edge $e$ is incident to two faces $f,f'$ the dual of $e$ in $P_w$ joins the vertices drawn in $f$ and $f'$.) For external permeable edges we need to be more precise. We consider a permeable edge $e$ incident to the root-face and a black-face $f$. The dual of $e$ joins the vertex drawn in $f$ to $v_i$ if $h_i$ is the last dangling head encountered before $e$ around the root-face, or to $v_0$ if no dangling head precedes $e$.
%The dual of $e$ join the vertex drawn in $f$ to $v_0$ if no dangling edge precedes $e$ during the tour of the root-face. It joins the vertex drawn in $f$ to $v_i$ if $h_i$ is the last dangling head encountered before $e$. 
Note that the partition-tree $P_w$ can be drawn in such a way that no edge of $P_w$ intersects another.  For instance, the partition-tree associated to $w=baa\B{b}bb\B{a}a$ is shown in Figure \ref{fig:partition-tree2}.  

Moreover the vertices of the partition-tree have an \emph{activity}. We call \emph{white} and \emph{black} the vertices of $P_w$ corresponding to  \emph{white} and \emph{black} faces of $\Ms_w$. The \emph{active} white vertices are $v_0,\ldots,v_k$. The \emph{active} black vertices are the vertices corresponding to rooting vertices of $M_w$ (see Subsection \ref{section:prefix-maps} where the notion of \emph{rooting vertex} is introduced). The other vertices are said to be \emph{inactive}. 

It remains to define the root of the partition-tree. Consider the first edge $e$ followed around the root-face of $\Ms_w$. It is a permeable edge. Its dual $e^*$ in $P_w$ joins the exterior vertex $v_0$ to the vertex drawn in the black face corresponding to the root-vertex of $M_w$. The root of $P_w$ is incident to $v_0$ and follows $e^*$ in counterclockwise direction around $v_0$. This root is indicated in  Figure \ref{fig:partition-tree2}. \\

\noindent Observe that, when $w=\Xi(M_T)$ is a parenthesis-shuffle, the map $M_w=\orient{M}^T$ has no dangling heads and the partition-tree $P_w$ is $\Upsilon\circ \varphi_1(\orientT{M})$.\\

We now relate the partition-tree $P_w$ to  the binary tree $\lambda_1(w)$.
\begin{prop}\label{thm:partition-tree-lambda1}
For all prefix-shuffle $w$, the partition-tree $P_w$ is equal to $\theta\circ\lambda_1(w)$ where $\lambda_1(w)$ is the binary tree defined in Definition \ref{def:lambdaP}  and $\theta$ is the mapping defined in Definition \ref{def:theta}.
\end{prop}
 
Proposition  \ref{thm:partition-tree-lambda1} implies that for any  parenthesis-shuffle $w=\Xi(M_T)$ we have $P_w=\theta\circ\lambda_1'(w)$. Given that  $P_w=\Upsilon \circ \varphi_1(\orientT{M})$, we obtain  $\varphi_1(\orientT{M})=\Theta \circ\lambda_1'(w)$.\\
%As observed above, if $w=\Xi(M_T)$ we have $P_w=\Upsilon \circ \varphi_1(\orientT{M})$. Therefore, Proposition \ref{thm:partition-tree-lambda1} implies that $\varphi_1(\orientT{M})=\Theta \circ\lambda_1'(w)$.\\

%\begin{cor}
%For all tree-rooted map $M_T$ of corresponding parenthesis-shuffle $w=\Xi(M_T)$, the non-crossing partitions  $\varphi_1(\orientT{M})$ and $\Theta\circ\lambda_1(w)$ are the same.
%\end{cor}

%\noindent \textbf{proof of the corolary:}
%Let $M_T$ be a tree-rooted map and $w=\Xi(M_T)$.  By Lemma \ref{thm:partition-tree-partition}, we have $P_w=\Upsilon \circ \varphi_1(\orientT{M})$. By Proposition \ref{thm:partition-tree-lambda1}, we have  $P_w=\theta\circ\lambda_1(w)$ hence $\Upsilon^{-1}(P_w)=\Theta \circ\lambda_1(w)$. Therefore, $\varphi_1(\orientT{M}) = \Theta\circ\lambda_1(w)$.
%\findem

The rest of this subsection is devoted to the proof of Proposition \ref{thm:partition-tree-lambda1}. We first describe a recursive construction of the partition-tree $P_w$. That is, we describe  how to obtain $P_{w\alpha}$ from $P_w$ when $\alpha$ is a letter in $\{a,\B{a},b,\B{b}\}$ (Lemma \ref{thm:partition-tree-evolution}). Then we describe a recursive construction of $\theta\circ\lambda_1(w)$ (Lemma \ref{thm:theta-rond-lambda}). We conclude the proof by induction on the length of $w$. \\

\subsubsection{Recursive construction of the partition-tree $P_w$}

The recursive description of the partition-tree requires to define an order on active vertices. Let $w$ be a prefix-shuffle and $M_w$ be the associated prefix-map.  The rooting vertices of $M_w$ can be compared by their order of  appearance around the spanning tree $T_w$ of $M_w$. The active black vertices inherit their order from the rooting vertices. The black vertex of $P_w$ corresponding to the root-vertex of $M_w$ is the first element for this order. We can also compare the dangling heads $h_1,\ldots,h_k$ of $M_w$ according to their order of appearance around  $T_w$. This order \emph{is the same} as the order of appearance around the root-face of $\Ms_w$.  Indeed, the order of appearance around the root-face of $\Ms_w$ is also the order of appearance around the root-face of $M_w$. Furthermore, the deletion of an edge of $M_w$ not in $T_w$ does not modify this order. By deleting all the edges not in $T_w$ we obtain the appearance order around $T_w$. The active white vertices inherit their order from the dangling heads. The exterior vertex $v_0$ is considered the first element. That is, $v_i$ \emph{precedes} $v_j$ for $0\leq i\leq j \leq k$.\\

%Let $w$ be a prefix-shuffle and $M_w$ the associated prefix-map. We call \emph{white} and \emph{black} the vertices of $P_w$ corresponding to  \emph{white} and \emph{black} faces of $\Ms_w$. The \emph{active} white vertices are $v_0,\ldots,v_k$. The \emph{active} black vertices are the vertices corresponding to rooting vertices of $M_w$. The other vertices are said \emph{inactive}. \\

%We need to order the active vertices. The rooting vertices of $M_w$ can be compared by their order of  appearance around the spanning tree $T_w$ of $M_w$. The active black vertices inherit their order from the rooting vertices. The black vertex of $P_w$ corresponding to the root-vertex of $M_w$ is the first element for this order. We can also compare the dangling heads $h_1,\ldots,h_k$ of $M_w$ according to their order of appearance around  $T_w$. This order \emph{is the same} as the order of appearance around the root-face of $\Ms_w$.  Indeed, the order of appearance around the root-face of $\Ms_w$ is also the order of appearance around the root-face of $M_w$. Furthermore, the deletion of an edge of $M_w$ not in $T_w$ does not modify this order. By deleting all the edges not in $T_w$ we obtain the appearance order around $T_w$. The active white vertices inherit their order from the dangling heads. The exterior vertex $v_0$ is considered the first element. That is, $v_i$ \emph{precedes} $v_j$ for $0\leq i\leq j \leq k$.\\

\begin{figure}[h!]
\begin{center}
\input{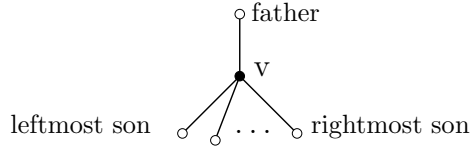}
\caption{A vertex and its leftmost and rightmost sons.} \label{fig:clockwise-son}
\end{center}
\end{figure}

Let $v$ be a vertex of a tree which is not a leaf. We call \emph{leftmost son} (resp. \emph{rightmost  son}) of $v$ the son following (resp. preceding) the father of $v$ (or the root) in counterclockwise direction around $v$ (see Figure \ref{fig:clockwise-son}).\\

We are now ready to describe the relation between the partition-tree $P_w$ and the partition-tree $P_{w\alpha}$ when $\alpha$ is a letter in $\{a,\B{a},b,\B{b}\}$.
\begin{lemma}\label{thm:partition-tree-evolution} 
The partition-tree $P_w$ is a tree. Moreover, 
\begin{itemize}
\item the partition-tree $P_{wa}$ is obtained from $P_w$ by adding a new leaf which becomes the last active black vertex. This leaf is  the leftmost son of the last active white vertex, 
\item the partition-tree $P_{wb}$ is obtained from $P_w$ by adding a new leaf which becomes the last active white vertex. This leaf is  the rightmost son of the last active black vertex,
\item the partition-tree $P_{w\B{a}}$ is obtained from $P_w$ by inactivating the last active black vertex,
\item the partition-tree $P_{w\B{b}}$ is obtained from $P_w$ by inactivating the last active white vertex.
\end{itemize}
\end{lemma}

\begin{figure}[h!]
\begin{center}
\input{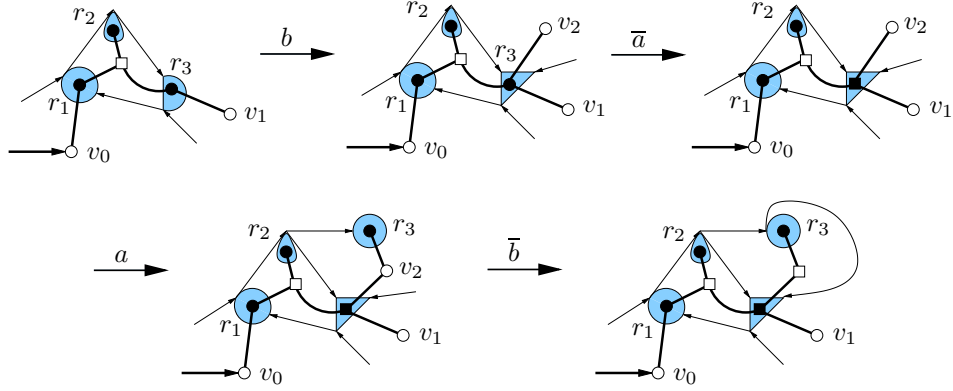}
\caption{Evolution of the partition-tree from $w=baa\B{b}b$ to $w=baa\B{b}bb\B{a}a\B{b}$.} \label{fig:partition-evolution}
\end{center}
\end{figure}

To illustrate this lemma we have represented the evolution of a partition-tree in Figure \ref{fig:partition-evolution}. Active vertices are represented by circles and inactive ones by squares. The white (resp. black) active vertices are denoted $v_0,v_1,\ldots$ (resp. $r_1,r_2,\ldots$).\\

\begin{figure}[h!]
\begin{center}
\input{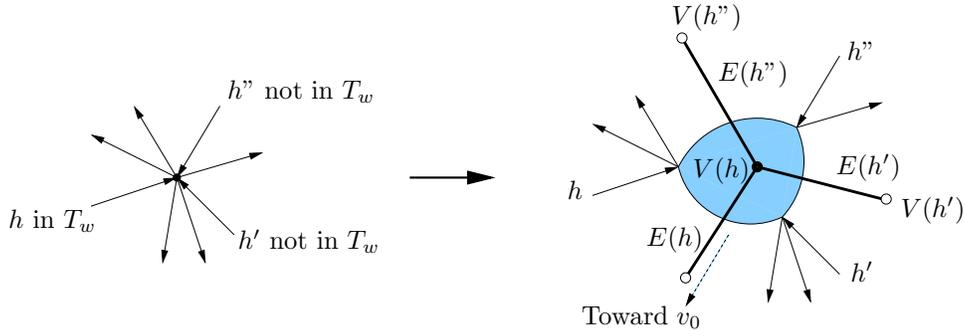}
\caption{Left: a typical vertex of the prefix map $M_w$ incident with three heads: $h$ in $T_w$ and $h'~,h''$ not in $T_w$. Right: the correspondence $E$ (resp.  $V$) between heads of $M_w$ and edges (resp. vertices) of $P_w$.  } \label{fig:EandV}
%\input{EandV.pstex_t}
%\caption{The edge $E(h)$ and the vertex $V(h)$ when $h$ is in $T_w$ (on the left) or not in $T_w$ (on the right).} \label{fig:EandV}
\end{center}
\end{figure}

Before we embark on the proof, we need to define a correspondence $E$ (resp. $V$) between the heads of $M_w$ and the edges (resp. vertices distinct from $v_0$) of $P_w$. The correspondences $E$ and $V$ are represented in Figure \ref{fig:EandV}.\\
Consider a head $h$ of $M_w$ and its end $v$ in $\Ms_w$. The edge following $h$ in counterclockwise direction around $v$ is a permeable edge. The dual of this edge in the partition-tree $P_w$ is denoted $E(h)$. The correspondence $E$ between heads of $M_w$ and edges of $P_w$ is one-to-one. The edge $E(h)$ is incident to a white and to a black vertex.  If $h$ is in the tree $T_w$ (in particular, if $h$ is the root), we define $V(h)$ as the black vertex incident to  $E(h)$. Else $V(h)$ is the white vertex incident to  $E(h)$.  The correspondence $V$ is a bijection between heads of $M_w$ and vertices of $P_w$ distinct from $v_0$. Indeed, black vertices of $P_w$ correspond to vertices of $M_w$ which are in one-to-one correspondence with heads in $T_w$, white vertices distinct from $v_0,\ldots,v_k$ correspond to faces of $M_w$ which  are in one-to-one correspondence with heads not in $T_w$ (a face $f$ is associated with the head we cross when we first enter $f$ during the tour of $T_w$), and the vertices $v_1,\ldots,v_k$ are  in one-to-one correspondence with the dangling heads $h_1,\ldots,h_k$. \\

%Consider a head $h$ of $M_w$ and its end $v$ in $\Ms_w$. The edge following $h$ in counterclockwise direction around $v$ is a permeable edge. The dual of this edge in the partition-tree $P_w$ is denoted $E(h)$. The correspondence $E$ between heads of $M_w$ and edges of $P_w$ is one-to-one. We now define a correspondence $V$ between heads of $M_w$ and vertices of $P_w$. If $h$ is in the tree $T_w$ (in particular, if $h$ is the root), the black vertex  $V(h)$ of $P_w$ corresponds to the vertex of $M_w$ incident to $h$. If $h$ is not in the tree $T_w$ and is not a dangling head, we consider the face $f$ of $M_w$ at its right. It is not the root-face of $M_w$. The white vertex $V(h)$ of $P_w$ is the vertex drawn in the face $f$. Lastly, if $h=h_i$ is a dangling head, the vertex $V(h)$ is $v_i$. This correspondence is a bijection between heads of $M_w$ and vertices of $P_w$ distinct from $v_0$. Indeed, black vertices of $P_w$ correspond to vertices of $M_w$ which are in one-to-one correspondence with heads in $T_w$, and white vertices distinct from $v_0,\ldots,v_k$ correspond to faces of $M_w$ which  are in one-to-one correspondence with heads not in $T_w$ (a face $f$ is associated with the head we cross when we first enter $f$ during the tour of $T_w$). \\

\dem We prove the lemma by induction on the length of $w$. If $w$ is the empty word, $P_w$ is a tree. Suppose now, by induction hypothesis, that $P_w$ is a tree. We first show  the following property: \emph{for any head $h$ of $M_w$, the edge $E(h)$ links $V(h)$ to its  father in $P_w$}. The mapping $V\circ E^{-1}$ is a bijection from the edges of $P_w$ to the vertices of $P_w$ distinct from its root-vertex $v_0$. Moreover an edge $e$ of $P_w$ is always incident to the vertex $V\circ E^{-1}(e)$ in $P_w$. Since $P_w$ is a tree,  the only possibility is that any edge $e$ of $P_w$ links the vertex $V\circ E^{-1}(e)$ to its father in $P_w$.

We are now ready to study separately  the different cases $\alpha=a,\B{a},b,\B{b}$. We use Lemma \ref{thm:prefix-map-evolution} and denote by $c$ the corner of $M_w$ at the left of the last rooting head of~$M_w$. 

\begin{figure}[ht!]
\begin{center}
\input{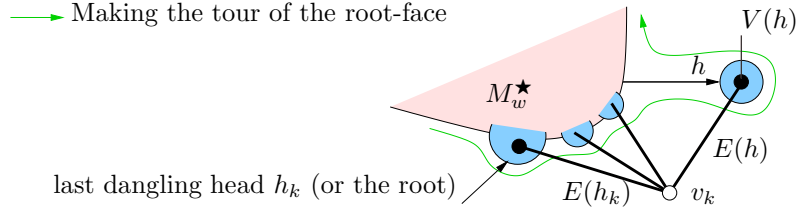}
\caption{The new vertex $V(h)$ is the leftmost son of $v_k$.} \label{fig:dem-prefix-mobile-a}\vspace{-.0cm}
\end{center}
\end{figure}

\ite Case $\alpha=a$. \\
\iten The prefix-map $M_{wa}$ is obtained from $M_w$ by adding a new edge $e$ in the corner $c$ oriented away from $c$. Let $h$ be the head of $e$ and $s$ its end. The vertex $s$ is the last rooting vertex in $M_{wa}$.
%The map $\Ms_{wa}$ is obtained  from $\Ms_w$ by adding the edge $e$ and the cell corresponding to $s$.  
The partition-tree $P_{wa}$ is obtained from $P_w$ by adding the edge $E(h)$ and the black vertex $V(h)$ to $P_w$ (see Figure \ref{fig:dem-prefix-mobile-a}). By definition, the vertex $V(h)$ is the last active black vertex in $P_{wa}$.\\
\iten By definition, the corner $c$ is situated after any dangling head around $T_w$. Hence, it is situated after any dangling head around the root-face of $\Ms_w$. Therefore, the edge $E(h)$ joins $V(h)$ to the last active white vertex $v_k$. Moreover, since $V(h)$ is only incident to $E(h)$ and $P_w$ is a tree, we check that $P_{wa}$ is a tree and $V(h)$ a leaf. \\
\iten It remains to show that $V(h)$ is the leftmost son of $v_k$. By definition, the permeable edges that have their dual incident to $v_k$ are situated between $h_k$ (or the root $h_0$ of $\Ms_w$ if $k=0$) and $c$ around the root-face of $\Ms_w$. The dual of the first of these permeable edge is $E(h_k)$ and the dual of the last of them is $E(h)$. If $k\neq 0$, we know that  $E(h_k)$ links $v_k=V(h_k)$ to its father in $P_w$. Therefore, $V(h)$ is the leftmost son of $v_k$. If $k=0$, we know  (by definition) that the root of $P_w$ follows $E(h_0)$ in counterclockwise direction around $v_0$. Therefore, $V(h)$ is the leftmost son of $v_0$.\\

\begin{figure}[ht!]
\begin{center}
\input{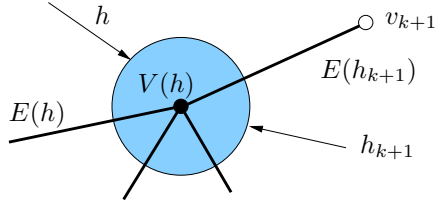}
\caption{The new vertex $v_{k+1}$ is the rightmost son of $V(h)$.} \label{fig:dem-prefix-mobile-b}\vspace{-.0cm}
\end{center}
\end{figure}

\ite  Case $\alpha=b$. \\
We denote by $h$ and $v$ the last rooting head and vertex.\\
\iten The prefix-map $M_{wb}$ is obtained from $M_w$ by adding a dangling head $h_{k+1}$ in the corner $c$.  It is the last dangling head of $M_{wb}$. The partition-tree $P_{wb}$ is obtained by adding  the vertex $v_{k+1}=V(h_{k+1})$ and the edge $E(h_{k+1})$ to $P_w$ (see Figure \ref{fig:dem-prefix-mobile-b}). By definition, $v_{k+1}$ is the last active white vertex of $P_{wb}$.\\
\iten The dangling head $h_{k+1}$ is incident to $v$ in $M_{wb}$. Hence,  the edge $E(h_{k+1})$ joins $v_{k+1}$ to the last active black vertex $V(h)$ of $P_{w}$. Moreover, since $v_{k+1}$ is only incident to $E(h_{k+1})$ and $P_w$ is a tree, $P_{wb}$ is a tree and $v_{k+1}$ a leaf.\\
%Since the corner $c$ is situated  at the left of the last rooting head $r$, the edge $E(h)$ joins $v_{k+1}$ to the last active black vertex $v$. Moreover, since $v_{k+1}$ is only incident to $E(h_{k+1})$ and $P_w$ is a tree, we check that $P_{wb}$ is a tree and $v_{k+1}$ a leaf. The leaf $v_{k+1}$ is the son of the last active black vertex $V(r)$.\\
\iten It remains to prove that $v_{k+1}$ is the rightmost son of $V(h)$. By definition, $E(h_{k+1})$ and $E(h)$ are respectively the dual of the permeable edges preceding and  following the head $h$  in counterclockwise direction around its end. Therefore, $E(h)$ follows $E(h_{k+1})$ in counterclockwise direction around $V(h)$. Given that $E(h)$ links $V(h)$ to its father, $v_{k+1}$ is the  rightmost son of $V(h)$.

\ite Case $\alpha=\B{a}$. \\
The prefix-map $M_{w\B{a}}$ is obtained from $M_w$ by inactivating the last active edge $e$.  Thus, $P_{w\B{a}}$ is obtained from $P_w$ by inactivating the last active black vertex.

\ite Case $\alpha=\B{b}$. \\
The prefix-map $M_{w\B{b}}$ is obtained from $M_w$ by adding a tail in the corner $c$ and connecting it to the last dangling head $h_k$. This creates a new face of $M_w$ (hence of $\Ms_w$) and lowers by one the number of dangling heads. The last active white vertex $v_k$ is  trapped in the new face of $M_{w\B{b}}$. Hence, $P_{w\B{b}}$ is obtained from $P_w$ by inactivating the last active black vertex $v_k$.
\findem

\subsubsection{Recursive construction of the tree $\theta \circ \lambda_1(w)$.}
We continue the proof of Proposition \ref{thm:partition-tree-lambda1}. We now describe the relation between the trees $\theta\circ \lambda_1(w)$ and $\theta\circ\lambda_1(w\alpha)$ when $\alpha$ is a letter in $\{a,\B{a},b,\B{b}\}$ (the mapping $\lambda_1$ is defined in Definition \ref{def:lambdaP}).\\

\begin{figure}[ht!]
\begin{center}
\input{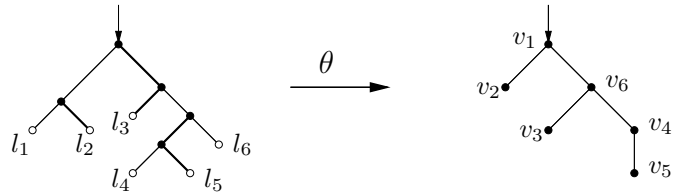}
\caption{Correspondence between leaves of $B$ and vertices of $\theta(B)$.} \label{fig:association-leaf-vertex}
\end{center}
\end{figure}

We first need to define a correspondence between the leaves of a binary tree $B$ and the vertices of the tree $\theta(B)$. 
An edge of $B$ is said \emph{left} (resp. \emph{right}) if it links a node to its left son (resp. right son). %A \emph{left path} (resp. \emph{right path}) is a maximal set of adjacent left edges (resp. right edges). 
We consider a leaf $l$ of $B$.  If $l$ is a left (resp. right) leaf, the path from $l$ to the root begins with a non-empty sequence of left (resp. right) edges.  By definition, only the last edge $e(l)$ of this sequence is branching except if $l$ is the first left leaf in which case no edge is branching. We associate the first left leaf of $B$ with the root-vertex of $\theta(B)$ and we associate any other leaf $l$ with the son of the branching edge $e(l)$ in $\theta(B)$. This correspondence is one-to-one. For instance, the leaves $l_1,\ldots,l_6$ of the binary tree $B$ in Figure \ref{fig:association-leaf-vertex} are associated with the vertices $v_1,\ldots,v_6$ of the tree $\theta(B)$.

Consider a prefix-shuffle $w$. In the binary tree $\lambda_1(w)$, leaves are either active or inactive. We say that a vertex of $\theta\circ \lambda_1(w)$ is left, right, active, inactive if the associated leaf of $\lambda_1(w)$ is so. Moreover, the leaves of the binary tree $\lambda_1(w)$ can be compared by their order of appearance around this tree. The vertices of $\theta\circ \lambda_1(w)$ inherit this order.  For instance, the root-vertex of  $\theta\circ \lambda_1(w)$ is the first active left vertex (recall that the first left leaf of $\lambda_1(w)$ is always active).\\

\noindent We are now ready to state the last lemma which is the counterpart of Lemma \ref{thm:partition-tree-evolution}.
\begin{lemma} \label{thm:theta-rond-lambda}
Let $\mathcal{T}$ be the tree $\theta\circ \lambda_1(w)$ and $\mathcal{T}_\alpha=\theta\circ \lambda_1(w\alpha)$ for $\alpha$ in $\{a,b,\B{a},\B{b}\}$.
\begin{itemize}
\item The tree $\mathcal{T}_a$ is obtained from $\mathcal{T}$ by adding a new leaf which becomes the first active right vertex. This leaf is  the leftmost son of the last active left vertex. 
\item The tree $\mathcal{T}_b$ is obtained from $\mathcal{T}$ by adding a new leaf which becomes the last active left vertex. This leaf is  the rightmost son of the first right vertex. 
\item The tree $\mathcal{T}_{\B{a}}$ is obtained from $\mathcal{T}$ by inactivating the first active right vertex.
\item The tree $\mathcal{T}_{\B{b}}$ is obtained from $\mathcal{T}$ by inactivating the last active left vertex.
\end{itemize}
\end{lemma}

\dem We study separately the four cases $\alpha=a,b,\B{a},\B{b}$.\\
\ite Case $\alpha=a$. By definition of the mapping $\lambda_1$ (Definition \ref{def:lambdaP}), the binary tree $\lambda_1(wa)$ is obtained from  $\lambda_1(w)$ by replacing the last active left leaf $l$ by a node with two leaves $l_l$ and $l_r$. 
%\centerline{\input{rule-P-a-dem-theta.pstex_t}}
The left leaf $l_l$ replaces $l$ as the last left leaf. The right leaf $l_r$ becomes the first right leaf. The edge from $l$ to $l_r$ is branching. The other branching edges are unchanged. Therefore,  $\mathcal{T}_a$ is obtained from $\mathcal{T}$ by adding a new leaf. This leaf is associated with $l_r$ hence becomes the first active right vertex. The father of this leaf was associated with $l$ in  $\mathcal{T}$ and is associated with $l_l$ in $\mathcal{T}_a$. Therefore, it was and remains the last active left vertex. It is easily seen that the new leaf becomes its leftmost son.  \\ 
\ite The case $\alpha=b$ is symmetric to the case $\alpha=a$. We do not detail it.\\
\ite  Case $\alpha=\B{a}$. The binary tree $\lambda_1(w\B{a})$ is obtained from  $\lambda_1(w)$ by inactivating the first active right leaf.  Therefore,  $\mathcal{T}_{\B{a}}$ is obtained from  $\mathcal{T}$ by inactivating the first active right vertex.\\
\ite  The case $\alpha=\B{b}$ is symmetric to the case $\alpha=\B{a}$.
\findem

\subsubsection{Recursive proof of Proposition \ref{thm:partition-tree-lambda1}.}
%We conclude the proof of Proposition \ref{thm:partition-tree-lambda1}. 
We want to show that, for any prefix-shuffle $w$, the partition-tree $P_w$ is the tree $\theta\circ \lambda_1(w)$. We show by induction the following more precise property: for any prefix-shuffle $w$, \\ 
\iten the partition-tree $P_w$ is equal to $\theta\circ \lambda_1(w)$ , \\
\iten the active and inactive vertices of $P_w$ and $\theta\circ \lambda_1(w)$ are the same,\\
\iten the white (resp. black) vertices of $P_w$ correspond to left (resp. right) vertices  of $\theta\circ \lambda_1(w)$,\\
%left (resp. right) vertices  of $\theta\circ \lambda_1(w)$ correspond to white vertices  (resp. black vertices) of $P_w$, \\
\iten the order on white (resp. black) vertices of $P_w$ is equal (resp. inverse) to the order on left (resp. right) vertices of $\theta\circ \lambda_1(w)$.\\

Suppose that $w$ is the empty word. The partition-tree  $P_w$ has one edge, an active white vertex which is its root-vertex and an active black vertex. Similarly, $\theta\circ\lambda_1(w)$ has one edge, an active left vertex which is its root-vertex and an active right vertex. Hence, we check that the property is true. 
%\begin{figure}[ht!]
%\begin{center}
%\input{P-epsilon.pstex_t}
%\caption{Partition-tree of the empty word.} 
%\end{center}
%\end{figure}
In view of Lemma \ref{thm:partition-tree-evolution} and Lemma  \ref{thm:theta-rond-lambda}, it is clear that the property is true by induction on the set of prefix-shuffles.
\findem

This concludes the proof of  Proposition \ref{thm:partition-tree-lambda1} and Theorem \ref{thm:equivalence}. \hspace{2.2cm} $\square$\\

\noindent \textbf{Acknowledgments:} I am deeply indebted to  Mireille Bousquet-Mélou for struggling with several versions of this paper and coming out with very helpful suggestions. This work has also benefited from fruitful discussions with Yvan Le Borgne, \'Eric Fusy and Gilles Schaeffer.

\bibliography{allref}
\bibliographystyle{plain}

\end{document}